\documentclass[english,11pt]{article}

\usepackage[german,french,english]{babel}
\usepackage[cp850]{inputenc}
\usepackage{latexsym,graphicx, fancybox}
\usepackage{graphicx}
\usepackage{amssymb, amsmath, amsfonts}
\usepackage{color}
\usepackage{latexsym}

\newcommand{\ve}{\varepsilon}  
  
\newcommand{\ds}{\displaystyle}  
\def\Q{{\rm I\hspace{-1.50ex}Q} }
\def\Q{\mathbb{Q}}
\def\R{\mathbb{R}}

\def\P{\mathbb{P}}
\def\E{\mathbb{E}}

\def\N{{\rm I\hspace{-0.50ex}N} } 
\def\N{\mathbb{N}}

\newtheorem{lem}{Lemma}[section]
\newtheorem{thm}{Theorem}[section]
\newtheorem{cor}{Corollary}[section]
\newtheorem{dfn}{Definition}
[section]
\newtheorem{prop}{Proposition}[section]

\begin{document}
\title{\bf Optimal Bregman quantization: existence and uniqueness of optimal quantizers revisited}

\author{ 
{\sc Guillaume Boutoille} \thanks{Laboratoire de Probabilit\'es, Statistique et Mod\'elisation, Sorbonne Universit\'e, case 158, 4, pl. Jussieu, F-75252 Paris Cedex 5, France  \& Fotonower, 30 rue Charlot, F-75003 Paris. E-mail: {\tt  guillaume.boutoille@orange.}}
\and   
{\sc  Gilles Pag\`es} \thanks{Laboratoire de Probabilit\'es, Statistique et Mod\'elisation,  Sorbonne Universit\'e. E-mail: {\tt  gilles.pages@upmc.fr}}~\thanks{This research
benefited from the support of the ``Chaire Risques Financiers'', Fondation du Risque.} }
	\maketitle
	
\begin{abstract} 
	In this paper we revisit the exsistence theorem for $L^r$-optimal quantization, $r\ge 2$,  with respect to a Bregman divergence: we establish the existence of optimal quantizaers under lighter assumptions onthe strictly convex function which generates the divergence, espcially in the quadratic case  ($r=2$). We then prove a uniqueness theorem  ``\`a la Trushkin'' in one dimension for strongly unimodal distributions and divergences gerated by strictly convex functions  whiose thire dervative is either stictly $\log$-convex or $\log$-concave.
\end{abstract}


%
%
%
%
%


\section{Introduction}In computer vision, labeling represents a major cost that we aim to reduce as much as possible. Clustering algorithms are efficient tools to partition similar data into clusters in order to organize data set. Browsing these clusters allows a better visualisation of the data set and an easier labeling. 
For that purpose, clustering is a fundamental field of "unsupervised" learning which aggregates procedures that have been widely studied across many disciplines. Most of the clustering methods consist in partitioning similar data into clusters usually ``represented'' or characterized by a {\em prototype} or {\em codeword} (either a typical true datum or a synthetic data). These codewords make up a {\em codebook}, such a codebook usually minimizes a loss function (quantization error). A widely used and studied clustering algorithm is the Euclidean $k$-means algorithm (\cite{macqueen67,jain88}). 

In Information Theory, Signal Processing and Probability Theory, this algorithm is also known as the Lloyd (or Lloyd~I) algorithm and goes back to the 1950s. At that time, it was a classified field of research in the Bell Labs, supported by the U.S. Navy (see the historical notes in the 1988 special issue of {\em IEEE on Information Theory} devoted to quantization). For probabilistic aspects, we refer to~\cite{GrafL2000}  who have proved most milestone theorems of the theory in a rigorous mathematical  way and, more recently,~\cite{LuPag23}. Since the 1970's (and even before), there  has been a  huge literature from Information Theory and automatic classification on the $k$-means algorithm  in connection with optimal quantization theory. Among them we can cite   the famous book by
~\cite{GershoG1991}.

In this paper we revisit the main existence theorem for the optimal quantizers of a distribution  $P$ supported by a (nonempty) open set $U\subset \R^d$  (in the sense that $P(U)=1$) when the loss function in the definition of the ($L^r$-mean) quantization error is  no longer a power of a norm on $\R^d$ but a power of a Bregman divergence induced  by a strictly convex  function $F: U \to \R$. Let us recall that the Bregman divergence $\phi_{_F}$  associated to  a strictly convex is defined on $U\times U\to \R_+$ by
\begin{equation}\label{eq:DefBregDiv}
\textcolor{black}{\forall\, (\xi,x)\!\in U^2, \quad \phi_{_F}(\xi,x)= \phi_{_F}(\xi,x):= F(\xi)-F(x)-\langle \nabla F(x)\,|\, \xi-x\rangle},
\end{equation}
where $\langle\cdot\,|\, \cdot \rangle$ denotes the  canonical inner product on $\R^d$. Indeed, the connection between these two frameworks is the case of the squared Euclidean norm which is its own Bregman divergence.  In~\cite{Banerjeeetal2005}, Banerjee and co-authors pointed out   a close relation between  Bregman divergences and $log$-likelihood of exponential families. This led them  to extend  the scope of $k$-means clustering algorithms by introducing  these Bregman divergences as new {\em  similarity measures} functions (a.k.a. {\em local loss functions}) in the core of the algorithm instead of a squared Euclidean norm. This work led to a line of research  about clustering based on Bregman divergences, see e.g.~\cite{Fischer2010},~\cite{HasFiscMou2021},~\cite{LiuBelNIPS2016}) among others.
An in-depth anamysis  of this extension of the $k$-means algorithm with Bregman divergnece as similarity measure is carried out in~\cite[Chapter 4]{bouto2024} and will soon be issued as a companion paper.
Nevertheless, we want first to  investigate  optimal vector quantization theory in the  same probabilistic spirit as that developed in  Graf  and Luschgy's book (\cite{GrafL2000}) or, more recently, in~\cite{LuPag23}, i.e. when the distribution $P$ to be quantized is a probability distribution rather than the empirical measure of a dataset (the data being modeled by a sample of i.i.d. random vectors in the statistical literature).

So we first establish an existence result for such $L^r$-optimal quantizers under some seemingly slightly weaker assumptions on the function $F$  than those  found in the (mathematically rigorous) literature (see e.g.~\cite{Fischer2010}).
We also investigate some properties of Bregman divergence based optimal quantizers known to be shared by regular optimal quantizers when the loss function is an Euclidean norm (see~\cite{GrafL2000, LuPag23}) like stationarity in the (pseudo-)quadratic case  and its $L^r$-variants when $r\neq2 $, but also the  seemingly less commonly known fact that the distribution $P$ under consideration assigns no mass to the boundaries of the Bregman--Voronoi cells of an optimal quantizer and other close properties. On the way, we provide existence result for the $(r,\phi_{_F})$-mean of $P$ under even less stringent assumptions on $F$. The $(r,\phi_{_F})$-mean of $P$ corresponds to a quantization of $P$ at level $n=1$. It is to be noticed that when $r=2$, the   $(2,\phi_{_F})$-mean of $P$ coincides with the standard mean of $P$ whatever $F$ is. 

The second main result, stated and established in Section~\ref{sec:BregTrushkin}, is that we extend for Bregman divergence based quantization the celebrated  uniqueness Trushkin theorem which states that, in one dimension and when $\phi_{_F}(\xi)= F(\xi)= |\xi|^2$,  if the distribution $P$ is strongly unimodal  (i.e.  is absolutely continuous with a $\log$-concave density on the real line), then stationary quantizers are unique. We prove that the uniqueness result remains true for our  extended Bregman setting under an additional assumption  on the function $F$, namely that its second derivative  \textcolor{black}{is either  $\log$-convex or $\log$-concave}. 

The 
paper
is organized as follows: Section~\ref{sec:Defs} is devoted to some short background on ``regular'' Optimal vector quantization  (by regular, we mean that  the loss function is a power of a norm) and the definitions of the mean quantization errors based on (powers of)  Bregman divergence $\phi_{_F}$. We also introduce the distortion which is the functional version of the ($r$th power) of the Bregman $(r,\phi_{_F})$-mean quantization error. Section~\ref{sec:optiBregquant},~\ref{sec:diffdistor} and \ref{sec:r-opti} are devoted to the existence theorem in the pseudo-quadratic case ($r=2$, see below) and, when $r> 2$, under appropriate integrability conditions on $P$ and light additional conditions on $F$ (on the boundary of $U$ and at infinity). We also establish the differentiability  of the distortion function, mostly  as a tool to establish the existence of Bregman $(r,\phi_{_F})$-optimal  quantizers when $r> 2$. Section~\ref{sec:BregTrushkin} is devoted to uniqueness of optimal quantizers at every level $n\ge 1$ when $P$ is strongly unimodal and $F$ is twice differentiable with $F''$ \textcolor{blue}{either  $\log$-convex or $\log$-concave} as well.

Our results are illustrated on several usual examples of Bregman divergences.

\section{Definitions and background on regular $L^r$-optimal quantization w.r.t. a norm}\label{sec:Defs}

\subsection{Short background on $L^r$-optimal vector quantization} 
 
Let  $\|.\|$ denote  a norm on $\R^d$. To any nonempty subset $\Gamma\subset \R^d$ we can associate Voronoi partitions of $\R^d$ induced by $\Gamma$.
\begin{dfn}[Voronoi partition] A Borel partition of $\R^d$ $(C_a(\Gamma))_{a\in \Gamma}$ is called a {\em Voronoi partition of $\R^d$ induced by $\Gamma$} if
\begin{equation}\label{eq:defVoro}
\forall\, a \!\in \Gamma,
C_a(\Gamma) \subset  \big\{\xi \!\in \R^d: \|\xi-a\| \le \|\xi-x\|,\; b\!\in \Gamma  \big\}.
\end{equation}
\end{dfn}

One checks that such a partition always  satisfies, for every $a\!\in \Gamma$,
\[
 \big\{\xi \!\in \R^d: \|\xi-a\| < \|\xi-b\|,\; b\!\in \Gamma\setminus \{a\}  \big\} \subset \mathring C_a(\Gamma) \subset  \bar C_a(\Gamma) \subset\big\{\xi \!\in \R^d: \|\xi-a\| \le \|\xi-b\|,\; b\!\in \Gamma  \big\}.
\]
Moreover, if $\|\cdot\|=|\cdot|_e$ is a Euclidean norm, the left and right inclusions hold as equality. Moreover the closed cells $C_a(\Gamma)$ are polyhedral convex sets i.e. the intersection of finitely many half-spaces (and $\mathring C_a(\Gamma)= \mathring{\bar C}_a(\Gamma) $). 

\medskip
Now let $P$ be a probability measure on $(\R^d,{\cal B}or(\R^d))\to \R^d$. Let $X:(\Omega, {\cal A}, \P)$ be an $\R^d$-valued random variable with distribution $P$. Let $r\!\in (0,+\infty)$. We can define the  quantization of $X$ induced by $\Gamma$ as follows.

Let us denote by ${\rm Proj}_{\Gamma}$ the {\em nearest neighbour projection} on $\Gamma$ induced by the partition $(C_a(\Gamma))_{a\in\Gamma}$ with respect to (w.r.t.)  the norm $\|\cdot\|$, defined  by 
\[
\forall\, x\!\in \R^d, \quad {\rm Proj}_{\Gamma}(x) =\sum_{a\in \Gamma}a\,\textbf{1}_{C_a(\Gamma)}(x).
\]
We make here an abuse of notation since such a projection not only depends on $\Gamma$  but also on the Voronoi partition  to which it is associated in a bijective way.   Note that for any Borel function $q:\R^d\to \Gamma$  
\[
\forall\, x\!\in \R^d, \quad \|x - q(x)\| \ge \| x-{\rm Proj}_{\Gamma}(x)\|
\]
with equality (on the whole space $\R^d$) if and only if  $q$ is a nearest neighbour projection  on $\Gamma$ in the sense that $\big(\{\xi: q(\xi)=a\}\big)_{a\in \Gamma}$ is a Voronoi partition of $\R^d$.

Then
\[
\widehat X^{\Gamma} := {\rm Proj}_{\Gamma}(X) =\sum_{a\in \Gamma}a\,\textbf{1}_{C_a(\Gamma)}(X)
\]
is  a {\em quantization of $X$} induced by the Voronoi partition $(C_a(\Gamma))_{a\in \Gamma}$ often called in short  a {\em quantization of $X$ induced by  $\Gamma$}. 

We can define likewise the distribution of $\widehat X^{\Gamma}$ characterized by its state space  $\Gamma$ and its  weights $\big(P(C_a(\Gamma))\big)_{a\in \Gamma}$ as the {\em quantization of the distribution $P$} of $X$ induced by the Voronoi partition $(C_a(\Gamma))_{a\in \Gamma}$ often called in short  a {\em quantization of the distribution $P$ induced by  $\Gamma$}, denoted  $\bar P$. 

Note that if $P$ has no atom then all quantizations of $X$ are $\P$-$a.s.$ equal whereas the quantization of $P$ is then unique. This is still true for distributions $P$ assigning no weights to hyperplanes  when the norm is Euclidean.

\begin{dfn}[$L^r$-quantization error induced by a quantizer] Let  $r>0$ and let $\|\cdot\|$ denote any norm on $\R^d$.  Assume that the distribution $P$ has a finite $r$th moment i.e. $\int_{\R^d}\|\xi\|^rP(\mathrm{d}\xi)<+\infty$. Let $\Gamma  \subset \R^d$ be a  {\em quantizer}. We define the $L^r$-mean quantization error of the distribution $P$ induced by $\Gamma$ by 
\begin{equation*}
    e_r(\Gamma,P,\|.\|) = \left[\int_{\mathbb{R^d}}\min_{a\in\Gamma} \|\xi-a\|^rP(\mathrm{d}\xi)\right]^{\frac{1}{r}}.
\end{equation*}
\end{dfn}

We can also define 
\begin{align*}
   e_r(\Gamma,X,\|.\|)   &:= \big\| \min_{a\in\Gamma} \|X-a\|\big\|_{L^r(\P)}\\
   &= \left[ \E\, \min_{a\in \Gamma}\|X-a\|^r\right]^{1/r}=   \left[\int_{\R^d}\min_{a\in\Gamma} \|\xi-a\|^rP(\mathrm{d}\xi)\right]^{1/r}\\
   &= e_r(\Gamma,P,\|.\|).
\end{align*}
Then, one straightforwardly checks that, for any  Voronoi partition, one has 
\[
\big\| \|X-\widehat X^{\Gamma}\|\big\|_{L^r(\P)}= \E\int \|{\rm Proj}_{\Gamma} (\xi)-\xi\|^rP(\mathrm{d}\xi)=  \sum_{a\in \Gamma}\int_{C_a(\Gamma)} \|\xi-a\|^rP(\mathrm{d}\xi) = e_r(\Gamma, P, \|\cdot\|)^r
\]
so that $e_r(\Gamma, P, \|\cdot\|)$ measures for the $L^r$-norm the  error induced when replacing a random variable $X$  by (any of) its (nearest neighbour) quantizations. 
This naturally leads to the following definition.
\begin{dfn}[$L^r$-optimal quantization error] 
The  $L^r$-optimal mean quantization error for $P$ at level $n\ge 1$ is defined by 
\begin{equation*}
    e_{r,n}(P,\|.\|) = \inf_{\Gamma: |\Gamma|\leq n}  e_r(\Gamma,P,\|.\|).
\end{equation*}
\end{dfn}
This in turn suggests the question of the existence of a solution to this minimization problem i.e. the question of the existence of an $L^r$-optimal quantizer of size (or {\em at level}) at most $n$ i.e. a quantizer $\Gamma$ which induces (through $\widehat X^{\Gamma}$) the best approximation of $X$ in $L^r(\P)$ by a random variable whose set of values has at most $n$ elements.

\medskip
\noindent {\bf Remark.} One easily checks (see e.g.~\cite{GrafL2000, LuPag23}) that 
\[
  e_{r,n}(P,\|.\|)= {\cal W}_r\big(P, {\cal P}_n\big),
\]
where ${\cal W}_r$ denotes the Wasserstein-Monge-Kantorovich distance between $P$ and the set of probability distributions on $\R^d$ whose support contains at most $n$ elements. In fact one also have for every $\Gamma\subset \R^d$
\[
  e_{r,n}(\Gamma,P,\|.\|)= {\cal W}_r\big(P, {\cal P}_{\Gamma}\big)
\]
where $ {\cal P}_{\Gamma}$ is the set of $\Gamma$-supported distributions.

This follows from the fact that, if $X\sim P$, then, for any random vector $\textcolor{black}{Y}:(\Omega, {\cal A}, \P)\longrightarrow \Gamma\subset \R^d$,
\begin{align*}
\big\| \|X-Y\| \big\|^r_{L^r(P)}&= \sum_{a\in \Gamma}\int_{\{Y=a\}}\|X-a\|^r\mathrm{d}\P\\
				 &\ge \sum_{a\in \Gamma}\int_{\{Y=a\}}\min_{b\in \Gamma}\|X-b\|^r\mathrm{d}\P\\
				 &= \E\, \min_{b\in \Gamma}\|X-b\|^r\mathrm{d}\P \\
				 &= \big\| \|X-\widehat X^{\Gamma}\| \big\|^r_{L^r(P)}= e_r(\Gamma,P)^r
\end{align*}

\subsection{$L^r$-optimal quantizers  for regular $L^r$-quantization: existence and properties}

%
This theorem is essentially proved in~\cite[Theorem~4.3, Chapter~1]{GrafL2000}. See also~\cite[Section~1.3]{LuPag23} for  claim~$(ii)$ of~$(b)$. 
\begin{thm}\label{thm:reg_existence} Let $r>0$ and let $\|\cdot\|$ denote any norm on  $\mathbb{R}^d$.  Assume $\int_{\mathbb{R}^d}\|\xi\|^{r}\P(\mathrm{d}\xi)< + \infty$.

\smallskip
\noindent $(a)$ For every $n\ge 1$, there exists at least one $L^r$-optimal quantizer $\Gamma^{(r,n)}$ of size at most $n$, solution to the minimization problem
\[
\inf_{\Gamma: |\Gamma|\leq n}  e_r(\Gamma,P,\|.\|)
\]
i.e. such that 
\[
e_r(\Gamma^{(r,n)},P,\|.\|)=  e_{r,n}(P,\|.\|).  
\]

\smallskip
\noindent $(b)$ If ${\rm card}\big({\rm supp}(P)\big) \ge n$, then any $L^r$-optimal quantizer $ \Gamma^{(r,n)}$ has full size $n$ and for every $a\!\in  \Gamma^{(r,n)}$, $P(C_a( \Gamma^{(r,n)}))>0$. Moreover, 

\smallskip
 $(i)$ if the norm $\|\cdot\|$ is smooth~(\footnote{that $\|\cdot\|$ is differentiable on $\R^d\setminus\{0\}$, like Euclidean norms or $\ell^p$-norms, $1<p<+\infty$.}) and strictly convex then {\em the boundary of the Voronoi partition is $P$-negligible}: for every $a\!\in \Gamma^{(r,n)}$, $P\big(\partial C_a(\Gamma^{(r,n)})\big)=0$

\smallskip
  $(ii)$ if   $r\ge 2$ and the norm $\|\cdot\|=|\cdot|_e$ is Euclidean then these quantizers are {\em $L^r$-stationarity}: for every $a\!\in \Gamma^{(r,n)}$,
\[
a = \frac{\int_{C_a(\Gamma^{(r,n)})}\xi|\xi-a|^{r-2}P(\mathrm{d}\xi)}{\int_{C_a(\Gamma^{(r,n)})}|\xi-a|^{r-2}P(\mathrm{d}\xi)}.
\]
\
\noindent $(c)$ When $d=1$  and $P=U([0,1])$, then for every $r>0$, the $L^r$-optimal quantizer is unique and 
\[
\Gamma^{(r,n)} = \Big\{\frac{2k-1}{2n},\; k=0,\ldots,n\Big\}.
\]
\end{thm}

 \noindent {\bf Remark.} In view of what follows with Bregman divergence it is interesting to inspect a variant of the above result as concerns the state space: \textcolor{black}{i}f the norm under consideration is an Euclidean norm denoted $|\cdot|_{e}$ and $P$ is supported by a (nonempty) closed convex set, say $C$, then it follows from the well-known fact that the resulting projection  $\pi_{C}$ on $C$  is $|\cdot|_{e}$-Lipschitz with Lipschitz coefficient $1$, that  for any quantizer $\Gamma\subset \R^d$
\[
e_r(\Gamma, P, |\cdot|_{e}) \textcolor{black}{\,\ge\,} e_r\big(\{\pi_{C}(a), \, a\!\in \Gamma\}, P, |\cdot|_{e}\big)
\]
since $|\xi-  \pi_{C}(a)|_{e}\le |\xi-a|_{e}$ for every $\xi\!\in C$ and $a\!\in\R^d$. Hence, one easily checks that  
\begin{align*}
e_{n,r}(P, |\cdot|_{e})^r &= \inf \big\{ e_r(\Gamma, P, |\cdot|_{e})^r : \Gamma\subset C, \; |\Gamma|\le n\big\}\\
& =  \inf \Big\{ \int \min_{b\in \Gamma}|\xi-b|_e^r P(\mathrm{d}\xi) : \Gamma\subset C, \; |\Gamma|\le n\Big\}.
\end{align*}
Now assume that $U=\mathring{C}\neq \varnothing$ and  that $P(U)=1$. 
Then, another argument (see among others~\cite{Pages1998},~\cite{GrafL2000} or~\cite{LuPag23}) based on the existence of optimal quantizers is to say that an $L^r$-optimal quantizer is always $U$-valued in that situation as a consequence of the support hyperplane theorem.

\section{Vector quantization with respect to a Bregman divergence: definitions and first properties} 
  \subsection{Bregman divergence induced by a strictly convex function $F$} First we introduce the notion of Bregman divergence induced by a strictly convex, continuously (Fr\'echet-)differentiable function $F$ defined on a nonempty convex open set $U$ of $\R^d$. This will be the class of similarity measures 
 that we will investigate in this paper.
  
\begin{dfn}[Bregman divergence associated to strictly convex function $F$]\label{def:Bregdiv}
Let $F:U\to \R^d$ be a continuously differentiable,  strictly convex function defined on a  {\em convex} $U\subset \R^d$. The convex $U$ will be assumed {\em open} in all what follows, except specific mention  (see e.g. Proposition~\ref{prop:existence}$(a)$). The $F$-Bregman divergence of $x\!\in U$ with respect to (w.r.t.) another point $\xi\!\in U$ is defined by 
\begin{equation}\label{eq:divBregmandef}
 \forall\, x\!\in U,\quad \phi_{_F}(\xi,x)= F(\xi)-F(x) -\langle \nabla F(x) | \xi- x \rangle\ge 0.
\end{equation}
Moreover $\phi_{_F}(\xi,x)=0$  if and only if $\xi=x$ since $F$ is strictly convex. 
\end{dfn}

When no ambiguity and to alleviate notations, we will often denote $\phi(\xi,x)$ instead of $\phi_{_F}(\xi,x)$.

\medskip It is important to notice that, if $F(x) = |x|^2 = x^* x$, $x\!\in \R^d$ 
one checks that $\phi_{_F}(\xi,x) =|\xi-x|^2$. This is also true for other Euclidean norms since if $F(x) =|x|_{_S}^2:= x^*Sx$ with $S\!\in {\cal S}^{++}(d,\R)$ (positive definite symmetric), then $\phi_{_F}(\xi,x)= |\xi-x|_{_S}^2$. This is the main connection between quantization w.r.t. to an $F$-Bregman divergence as a similarity measure (or local loss function) and that w.r.t. to a norm (or their $r$-powers).

\medskip Various examples  of Bregman divergence  of interest for applications are presented in Section~\ref{subsec:existr=2} right after the first Theorem~\ref{thm:Existence}  of existence of optimal $\phi_{_F}$-quantizers at level $n\ge 1$ of  a distribution $P$.

\medskip 
\noindent {\bf Remarks.} $\bullet$ 
When $F$ is twice (Fr\'echet-)differentiable on $U$ then a second order Taylor expansion with integral remainder yields the alternative formulas that is also useful fro many applications
\begin{equation*}\label{eq:bregmadivFC2}
   \phi_{_F}(\xi,x)= \int_0^1 (1-t)\nabla^2 F\big(x+t(\xi-x)\big).(\xi-x)^{\otimes 2} \mathrm{d}t=  \int_0^1s\nabla^2 F\big(\xi+s(x-\xi)\big).(x-\xi)^{\otimes 2} \mathrm{d}s.
\end{equation*}

\noindent $\bullet$ Generally Bregman divergences are not symmetric and do not satisfy the triangle inequality, which makes them significantly different from a distance since it only satisfies {\em a priori} $\phi_{_F}(\xi,x)= 0$ if and only of $\xi=x$.
    That can be paradoxically summed by: for $\xi$, $x$, $y\! \in U$, the  following usually holds
     \smallskip
     \begin{itemize}
 	 \item $\phi_{_F}(\xi,x) \neq \phi_{_F}(x,\xi)$.
        \item $\phi_{_F}(\xi,x) \nleq \phi_{_F}(\xi,x) + \phi_{_F}(x,y)$,
 \end{itemize}
 
 \smallskip
  However Bregman divergences share some somewhat hidden properties that (partially) compensate this problem.  Enough to be quite  convincing as  similarity measure to define and develop a optimal quantization/clustering theory.
  Throughout this manuscript, the function $F$ will denote a strictly convex ${\cal C}^1$-function defined on a (nonempty) convex open set $U$ having values in $\R$, the  divergence $\phi_{_F}$  defined by~\eqref{eq:divBregmandef} will denote the resulting $F$-Bregman divergence.
  
 \begin{prop}\label{prop:bgredivprop1}  Let $F$ and $G$ denote two ${\cal C}^1$ strictly convex functions. 
 
 \smallskip
 \noindent $(a)$ Then for every $\lambda$,  $\mu \!\in \R_+$, $\lambda+\mu \neq 0$,  $\lambda F+\mu G$ is still ${\cal C}^1$ and strictly convex. Furthermore
 \[
 \phi_{_{\lambda F+\mu G}} = \lambda \phi_{_F} +\mu \phi_{_G}.
 \]

\noindent $(b)$ One has 
\[
\phi_{_F}= \phi_{_G} \Longleftrightarrow \exists\, \alpha \! \in \R^d, \; \beta \!\in \R \;\mbox{ such that }\;\forall\,\xi
\!\in \R^d,\; F(\xi)-G(\xi) = \langle \alpha | \xi\rangle +\beta.
\]
\end{prop}

\noindent {\em Proof.} $(a)$ is obvious. 

\smallskip
\noindent $(b)$ The direct sense follows from the fact that formally $\phi_{\langle\alpha|\cdot\rangle+\beta}= 0$ on $U$. Conversely, if $\phi_{_F}= \phi_{_G}$ then, any $x_0\!\in \R^d$ being fixed,  the announced identity holds for every $\xi \!\in \R^d$ with $\alpha = \nabla (F-G)(x_0)$ and $\beta = (F-G)(x_0)+ \langle \nabla (F-G)(x_0)\,|\, x_0\rangle$.
\hfill$\Box$

\smallskip
For some Lipschitz/strongly convex functions $F$ it is possible to compare the Bregman divergence $\phi_{_F}$ to the Euclidean norm  which is its own Bregman divergence as seen above.
\begin{prop}[$F$-Bregman divergence compared to Euclidean norm]\label{property:2}
    $(a)$ Assume that $F:U \to \R$ has a Lipschitz continuous gradient $\nabla F$ on $U$. Then, for every  $\xi,\,x\!\in U$,
    \begin{equation}\label{eq:ComparEuclide}
        0 \leq \phi_{_F}(\xi,x) \leq \tfrac 12  [\nabla F]_{\rm Lip}|\xi- x|^2.
    \end{equation}
    
    \noindent $(b)$ If $F$ is $\alpha$-convex for some $\alpha >0$ in the sense that
    \[
   \forall\, \xi, x\!\in U, \quad  \langle \nabla F(\xi ) - \nabla F(x) | \xi- x \rangle \ge \alpha |\xi- x|^2
    \]
    then 
    \[
     \phi_{_F}(\xi,x) \ge  \frac{\alpha}{2}  |\xi- x|^2.
     \]
     \end{prop}
     
\noindent {\em Proof.} $(a)$ Let $\xi, x\!\in U$. One has using Cauchy-Schwartz inequality in the second line 
    \begin{align*}
        \phi_{_F}(\xi,x) &= F(\xi) - F(x) -\langle \nabla F(x) | \xi-x \rangle \\
       & =\int_0^1  \langle \nabla F((1-u)\xi +u x) -\nabla F(x)\,|\, \xi-x\rangle du\\
    & \le \int_0^1  \big|\nabla F((1-u)\xi +u x) -\nabla F(x)\big|| \xi-x| du\\
    & \le  \int_0^1  (1-u)du [\nabla F]_{\rm Lip} | \xi-x|^2  =  \tfrac 12 [\nabla F]_{\rm Lip}|\xi-x |^2.
    \end{align*}
    
    \noindent $(b)$ is an easy consequence of the fact that the function
    \[
    g(t)  = F(t\,\xi+(1-t)x)-F(x) -t\langle \nabla F(x)\,|\, \xi-x\rangle -\frac{\alpha}{2}t^2|\xi-x|^2
    \]
     has a non-negative derivative on $(0,1]$ since 
     \begin{align*}
     g'(t) &= \langle \nabla F(t\xi+(1-t)x)-\nabla F(x)\,|\,\xi-x\rangle -\alpha t|\xi-x|^2\\
     & =\frac 1t \Big(\big\langle \nabla F(t\xi+(1-t)x)-\nabla F(x)\,|\,t(\xi-x)\big\rangle -\alpha \big|t(\xi-x)\big|^2\Big)\ge 0
     \end{align*}
     which implies  $g(1)\ge g(0)=0$. 
     \hfill$\Box$
     
      \begin{prop}[Integrability]  Let $r>0$. If  $\mathbb{E}\big(|F(X)| \vee |X|\big)^{\frac r2}  < + \infty $, then for every $x\!\in U$, $\mathbb{E}\ \phi_{_F}(X,x)^{\frac r2} < +\infty$.
\end{prop}

\noindent {\em Proof.} Combining Cauchy-Schwartz inequality and the elementary inequality $|a+b|^{\frac r2} \le 2^{(\frac r2-1)^+}\big(|a|^{\frac r2} +|b|^{\frac r2} \big)$, yields 
    \begin{align*}
      \hskip 3cm  \mathbb{E}\,|\phi_{_F}(\xi,x)|^{\frac r2 }&= \int_U |F(\xi) -F(x) - \langle\nabla F(x)| \xi-x\rangle | ^{\frac r2}P(\mathrm{d}\xi) \\
        & \le 2^{(\frac r2-1)^+}(1+|\nabla F(x)| )^{\frac r2} \int_U \big(|F(\xi) | \vee  | \xi| \big)^{\frac r2}P(\mathrm{d} \xi) \\
        & \hskip 3cm  + |F(x)|^{\frac r2} +|\langle \nabla F(x) | x\rangle |^{\frac r2} <+ \infty.\hskip 3cm \Box
    \end{align*}   
  \subsection{Quantization with respect to an $F$-Bregman divergence}  
    
    The idea of Bregman quantization is to replace the norm  which plays the role of a similarity measure in regular vector quantization theory  by the   Bregman divergence $\phi_{_F}$ of a continuously  differentiable  strictly convex function $F$ as defined  in~\eqref{eq:bregmadivFC2}. The  function $F$ being defined on a nonempty  convex  open set $U\subset \R^d$, we will consider in what follows $U$-valued quantizers in our definitions of the quantization  error. This choice, although natural,  may appear somewhat arbitrary when the support of $P$ is strictly included in $U$ or, even worst, if $U$ contains the (closed)  convex hull of the support of $P$. Except in $1$-dimension (see Section~\ref{sec:BregTrushkin}), it seems not clear that it  has no influence on the  quantization error, compared to another natural choice which could be the convex hull  of the support of the distribution $P$ (or preferably its interior if it has $P$-probability $1$).

\begin{dfn}[$(r,\phi_{_F})$-Quantization Error]\label{def:BregQerrorchap2} 
Let $r\!\in (0, +\infty)$.  Let $P$ be a probability distribution supported by $U$ and let $X:(\Omega,{\cal A}, \P)\to \R^d$ be a $P$-distributed random vector. Assume
\begin{equation}\label{eq:IntegrabiliteBregquant}
\mathbb{E}\big(|F(X)| \vee |X|\big)^{\frac r2}  < + \infty.
\end{equation}
\noindent $(a)$  Let $\Gamma\subset U$ be  a finite subset of $U$ (also called {\em quantizer}). The $(r,\phi_{_F})$-mean quantization error for the Bregman divergences $\phi_{_F}$ of the distribution $P$ induced by $\Gamma$ is defined by 
\begin{equation}\label{eq:meanBregQerrorchap2}
    e_{r,n}(\Gamma, P,\phi_{_F}) =  \left[\int_U \min_{a \in \Gamma}  \phi_{_F}(\xi,a)^{\frac r2} P(\mathrm{d}\xi) \right]^{\frac{1}{r}}= \left[ \E \, \min_{a\in \Gamma}\phi_{_F}(X,a)^{\frac r2}\right]^{\frac 1r}<+\infty.
\end{equation}

\smallskip
\noindent $(b)$ The optimal $(r,\phi_{_F})$-mean quantization error for Bregman divergences at level $n\!\in \N$ is defined by
\begin{equation}\label{eq:OptimeanBregQerrorchap2}
    e_{r,n}(P,\phi_{_F}) = \inf_{ |\Gamma|\leq n}  e_{r,n}(\Gamma, P,\phi_{_F}).
\end{equation}
\end{dfn}

We will also use   the notation    $e_{r}(\Gamma, X,\phi_{_F})$ and  $e_{r,n}(X,\phi_{_F})$ to denote the above quantities, depending on the context; keeping in mind that these quantizers only depend on the distribution $P$, not on the random vector with distribution $P$ under consideration.

\bigskip
\noindent {\bf Comment.} The terminology  $(r,\phi_{_F})$-mean quantization error has been assigned to this error modulus to be consistent  with the usual terminology $L^r$-quantization when $F(x) = |x|^2$ (Euclidean norm). But this choice is somewhat arbitrary.

The quantities $e_{r}(\Gamma, X,\phi_{_F})$ can be defined regardless of the integrability condition~\eqref{eq:IntegrabiliteBregquant}, but when equal to $+\infty$ it is of no interest.


\section{Switch to the distortion function}\label{sec:distor}
\subsection{Definition of   $(r,\phi_{_F})$-distortion functions} 
%
Let $P$ be a probability distribution supported by the non empty open convex set $U$ i.e. such that $P(U)=1$, let $X$ be a random vector with distribution $P$ and let $r>0$ be such that
\begin{equation}\label{eq:Integrabilite}
\mathbb{E}\,\big(|X| \vee |F(X)| \big)^{\frac  r2}  =  \int_U \big(|\xi|\vee | F(\xi)| \big)^{\frac  r2}P(\mathrm{d}\xi) <+\infty.
\end{equation}

 We introduce, for every {\it level } $n\ge 1$, the   $(r,\phi_{_F})$ -distortion function $G^{F,P}_{r,n}:U^n\to \R_+$ of  the distribution $P$ and the Bregman divergence $\phi_{_F}$. The function $G_{r,n}= G_{r,n}^{F,P}$ is defined on $U^n$ by 
\begin{equation}\label{eq:n-distortion}
\forall \, x=(x_1,\ldots, x_n)\!\in U^n,\quad G_{r,n}(x_1,\ldots,x_n) = \int_U \min_{i=1,\ldots,n}\phi_{_F}(\xi,x_i)^{\frac  r2}P(\mathrm{d}\xi).
\end{equation}

It follows from~\eqref{eq:divBregmandef} that for every $x \!\in U$, there exists $C_x>0$ such that  for every $\xi \!\in U$
\begin{equation}\label{eq:boundphi1}
0\le \phi_{_F}(\xi,x) \le C_x\big(1+|F(\xi)|\vee|\xi| \big),
\end{equation}
where $C_x$ continuously depends on $|F(x)|\vee |x|\vee|\nabla F(x)|$ so that the distortion function $G_{r,n}$ is well defined as an $\R_+$-valued function defined on $U^n$. This assumption will be in force in all what follows.

\bigskip
\noindent {\bf Remark.} When $F= |\cdot|^2$, then $\phi_{_F}(\xi,x)= |\xi-x|^2$ and $G_{r,n}$ is simply the regular $L^r$ distortion function with respect to the Euclidean norm. This is the reason for choosing $\frac r2$ as an exponent rather than simply $r$.

\bigskip
For notational convenience, when $r=2$ we will drop   the dependence in $r$  and write $G_n$ instead of $G_{r,n}$. 

\smallskip It is also  convenient to rely on a more probabilistic formulation based on random variables. Thus,  if $X: (\Omega,{\cal A}, \P)\to \R^d$ with distribution $P$, then~\eqref{eq:Integrabilite} reads equivalently
\begin{equation}\label{eq:Integrabilite-bis}
\E\,\big(|X| \vee |F(X) |\big)^{\frac  r2}   <+\infty
\end{equation}
and the  $(r,\phi_{_F})$-distortion function  reads
\begin{equation}\label{eq:n-distortion2}
\forall\, (x_1,\ldots,x_n) \!\in U^n, \quad G_{r,n}(x_1,\ldots,x_n)  = \E\, \min_{i=1,\ldots,n} \phi_{_F}(X,x_i)^{\frac  r2}.
\end{equation}
We will use both formulations.

 \subsection{Properties of $F$-Bregman divergences and $(r,\phi_{_F})$-distortion functions}
 
 We will extensively use the following properties satisfied by $P$, $\phi_{_F}$ and $G_{r,n}$.
 
 \bigskip 
 \noindent {\bf P1.} If $P$ has a first moment. Then
 $$
 \E\, X =\int_U\xi P(\mathrm{d}\xi) \!\in U
 $$ 
since  $U$ is convex.  

\smallskip
In fact if $\E\, X \!\in \partial U= \bar U \setminus U$, then  the supporting hyperplane theorem applied at $x=\E\, X\!\in \partial U$ to the convex set $U$ yields a contradiction. To be more precise, there  exists a vector $u^{\perp}\!\in \R^d$, $|u^{\perp}|=1$ such that $U\subset \{\xi: \langle\xi-x\,|\,u^{\perp} \rangle <0\}$ so that integrating w.r.t. $P$ implies $\big\langle \int_U\xi P(\mathrm{d}\xi)-x\,|\,u^{\perp}\big\rangle <0$ which is impossible.

 \medskip 
 \noindent {\bf P2.} {\em A bound of interest.} For every $\xi,\,x\!\in U$, 
\begin{equation}\label{eq:n-distortion3}
0\le \phi_{_F}(\xi,x) \le \langle \nabla F(x)-\nabla F(\xi)\,|\, x-\xi\rangle.
\end{equation}
This follows from the convexity inequality $F(x) \ge F(\xi)-\langle \nabla F(\xi)\,|\, x-\xi\rangle$ satisfied by the function $F$. 

\smallskip This provides a criterion for the  finiteness of the  $ (r,\phi_{_F})$-distortion function $G_{r,n}$, namely
\[
\E\, (|\nabla F(X)||X| )^{r/2}= \int_U (|\nabla F(\xi)||\xi|)^{r/2}P(\mathrm{d}\xi)<+\infty, 
\]
since
$$
|F(\xi) |\le | F(\xi_0) |+ |\langle \nabla F(\xi)\,|\, \xi-\xi_0\rangle |\le  | F(\xi_0) |+ |\nabla F(\xi)| |\xi-\xi_0|
$$ 
for any $\xi, \,\xi_0\!\in U$.

 \medskip 
 \noindent {\bf P3.} {\em A useful identity.} The Bregman divergence is not a distance in general since it is usually not symmetric and, even more penalizing, it does not satisfy a triangle inequality. However the following elementary identity can be sometimes called upon  to ``make the job''. (We leave its easy  proof to the reader). For every $u,\,v, \, w\!\in U$, 
\begin{equation}\label{eq:phiabc}
\phi_{_F}(u,v)+ \phi_{_F}(v,w) = \phi_{_F}(u,w) + \langle \nabla F(w)-\nabla F(v)\,|\,u-v \rangle.
\end{equation}
which can also be rewritten
\begin{equation}\label{eq:phiabc2}
 \phi_{_F}(u,w)  = \phi_{_F}(u,v)+ \phi_{_F}(v,w) - \langle \nabla F(w)-\nabla F(v)\,|\,w-v \rangle+ \langle \nabla F(w)-\nabla F(v)\,|\,w-u \rangle.
 \end{equation}

  \noindent {\bf P4.} {\em Continuity of $G_{r,n}$ on $U^n$.} It is clear that for every $\xi\!\in U$,  $U^n \owns x=(x_1,\ldots,x_n)\mapsto \displaystyle \min_i\phi_{_F}(\xi,x_i)$ is continuous  on $U^n$ since $F$ is ${\cal C}^1$.  

Moreover on a compact neighbourhood   $K\subset U^n$ of $(x_1,\ldots, x_n)$, one has
\[
\forall\, y=(y_1, \ldots,y_n) \!\in K, \quad      \min_i\phi_{_F}(\xi,x_i)^{\frac  r2}\le\phi_{_F}(\xi,x_1)^{\frac  r2}\le C_K(|F(\xi)|\vee |\xi|)^{\frac  r2}
\]  
owing to~\eqref{eq:boundphi1} since $F$ and $\nabla F$ are bounded on $K$ by continuity. The  continuity of $G_{r,n}$ on $U^n$ follows by Lebesgue's dominated  continuity  theorem.

  \medskip
\noindent {\bf Remark.}   The same reasoning  shows that, if both $F$ and $\nabla F$ admit a continuous extension on $\bar U$, then,  $G_{r,n}$ can be defined for any distribution $P$ supported by $\bar U$ and continuously extended on $\bar U^n$ as well.  However, we will consider in what follows only distributions $P$ supported by $U$.

  \medskip
  \noindent {\bf P5.} {\em Lower semi-continuous extensions of $\phi_{_F}(\xi,\cdot)$ and $G_{r,n}$.}  
 
\noindent  $\rhd$ {\em General case.} For technical purposes we will extensively need to extend the functions $\phi_{_F}(\xi, \cdot)$, $\xi \!\in U$, over the Alexandroff~--~or one-point~--~compactification $\bar U^{\widehat{\R^d}}= \R^d\cup\{\infty\}$ of $U$ where $U$ is viewed  as an open subset of the Alexandroff compactification $\widehat{\R^d}$ of $\R^d$, equipped with its usual (metric) topology. Using this trick goes back to~\cite{Sabin-Gray} and is used in~\cite{Fischer2010} as well. If $U$ is bounded, then $\bar U^{\widehat{\R^d}}= \bar U$ and $\partial \bar U^{\widehat{\R^d}}= \partial U$  \textcolor{black}{otherwise}  i.e. $\bar U^{\widehat{\R^d}}= \bar U \cup\{\infty\}$. If we denote by
 $$
 \widehat{\partial} \,U :=\partial \bar U^{\widehat{\R^d}}
 $$ 
 the boundary of $U$ in $\widehat{\R^d}$, then $\widehat{\partial} \,U= \partial U$ if $U$ is  bounded and $\widehat{\partial} \,U= \partial U\cup\{\infty\} $ if $U$ is unbounded. 
 
 One extends the Bregman  divergences functions $\phi_{_F}(\xi, \cdot)$ for every $\xi\!\in U$ by their lower-semi-continuous  (l.s.c.) envelopes, still denoted $\phi_{_F}(\xi, \cdot)$,  defined on the boundary of $U$ by
 \begin{align}
\nonumber  \forall\,\bar x \!\in \widehat{\partial}  U,\;   \phi_{_F}(\xi, \bar x)  &:= \liminf_{x\to \bar x, \, x\in U}\phi_{_F}(\xi,x)\\
\label{eq:lscextphi} & := \inf\!\big\{\lim_n \phi_{_F}(\xi,x_n), \, x_n\to x,\,x_n \!\in U,\, \phi_{_F}(\xi,x_n) \mbox{ converges in } \bar \R_+\big\}. 
\end{align}
 
Note that this infimum is in fact a minimum and that one may consider  {\em a posteriori} $\bar U^{\widehat{\R^d}}$-valued sequences $(x_n)$ in the definition. Thus extended, $\phi_{_F}(\xi,\cdot)$ is l.s.c. for every $\xi\!\in U$.

\smallskip Then one extends accordingly the  $(r,P,F$)-distortion function $G_{r,n}$   to  $(\bar U^{\widehat{\R^d}})^n$ formally by~\eqref{eq:n-distortion} or~\eqref{eq:n-distortion2} using the above l.s.c. extension~\eqref{eq:lscextphi}   of the functions $\phi_{_F}(\xi,\cdot)$, namely
\begin{equation}\label{eq:lscGrn}
\forall\, \xi \!\in U,\;  \forall\, (x_1, \ldots,x_n) \!\in (\bar U^{\widehat{\R^d}})^n, \quad G_{r,n}(x_1, \ldots,x_n) = \E\, \min_{i=1,\ldots,n} \phi_{_F}(X,x_i)^{\frac r2}.
\end{equation} 

One checks that, as a straightforward   consequence of  Fatou's lemma,  such an extension of $G_{r,n}$  is  in turn l.s.c. on $(\bar U^{\widehat{\R^d}})^n$.

\smallskip $\rhd$ {\em The unidimensional case ($d=1$)}. When $d=1$ there is an alternative to the one-point Alexandroff's compactification of $\R$ that takes advantage of the natural order on $\R$. It consists in topologically identifying $\R$ to $(-1,1)$ by an appropriate increasing homeomorphism so that $\bar \R= [-\infty, \infty]$ is homeomorphic to $[-1,1]$. This compactification preserves the natural order on $\R$. The closure $\bar U ^{\bar \R}$ of  an unbounded convex interval $U\subset\R$  is obtained by adding to its closure $\bar U$ in $\R$, $\pm\infty$ depending on  possible infinite endpoints.  Then, one can extends the function $\phi_{_F}(\xi,\cdot)$, $\xi\!\in U$ by l.s.c. envelope  at every point  of the  boundary $\bar \partial\, U$ of $U$ in $\bar \R$.

  \medskip
  \noindent {\bf P6.} {\em Marginal convexity and monotonicity of Bregman divergence.}
 
\smallskip
 \noindent $\rhd$ {\em Convexity of $\phi_{_F}(\cdot,x)$.} First notice that for any fixed $x\!\in U$, $\xi \mapsto \phi_{_F}(\xi,x)$ is obviously convex as a sum of a convex and an affine  functions. Note that, by contrast, the function $x \mapsto \phi_{_F}(\xi, x)$ is not convex in general as illustrated by Figure~\ref{fig:phixinotconvex}.  
 
\medskip
 \noindent $\rhd$ {\em Monotonicity of  $\phi_{_F}(\xi,\cdot)$.} The following lemma deals with the monotonicity properties of $x \mapsto \phi_{_F}(\xi,x)$ for a fixed $\xi \!\in U$. For every point $\xi \!\in U$ and every vector $v\!\in S(0,1)$ (unit Euclidean sphere), we define
 \[
 t^-_{\xi,v}:= \inf\{t<0: \xi+tv \!\in U \}\!\in [-\infty,0)\quad\mbox{ and }\quad  t^+_{\xi,v}:= \sup\{t>0: \xi+tv \!\in U \}\!\in (0, +\infty].
 \]
 The signs of $t^\pm_{\xi,v}$ are clear since $U$ is open. Furthermore, for every $t\!\in (-t^-_{\xi,v},t^+_{\xi,v})$, $\xi+tv \!\in U$. 
 \begin{lem}[Monotonicity] \label{lem:monotonie} $(a)$ For every $\xi\!\in U$, $v\!\in S(0,1)$, the mapping $(t^-_{\xi,v}, t^+_{\xi,v}) \owns t \mapsto \phi_{_F}(\xi, \xi+t v)$ is decreasing on $(t^-_{\xi,v},0]$ and increasing on $[0,t^+_{\xi,v})$.

\noindent $(b)$  For every $\bar x \!\in \widehat \partial\, U$, $\xi\!\in U$, $\phi_{_F}(\xi, \bar x):= \displaystyle \liminf_{x\to\bar x, \,x\in U} \phi_{_F}(\xi,x) >0$. 

In particular, when $U$ is unbounded (so that $\infty\!\in \hat \partial U$)
\begin{equation}\label{eq:liminf>0}
\forall\, \xi\!\in U,\quad \liminf_{|x|\to \infty, \,x\in U}\phi_{_F}(\xi,x)>0.
\end{equation}

\noindent $(c) $ When $d=1$,  Claim $(b)$ is also true for the boundary $\bar \partial\, U$ of $U$ in $\bar \R= [-\infty, +\infty]$.
\end{lem}

\noindent {\em Proof.} $(a)$ First note that  the function $f(t)= \phi_{_F}(\xi, \xi +t\,u)$ is well defined on $(t^-_{\xi,u},t^+_{\xi,u})$. Let $t'>t>0$. 
 Elementary computations show that, for $t,\, t' \!\in [0, t^+_{\xi,u})$,  $t'>t$, 
 \small
 \begin{align*}
 \phi_{_F}(\xi, \xi+t'u)- \phi_{_F}(\xi, \xi+tu)&= F(\xi+t\, u) -F(\xi+t'u) + t'\langle \nabla F(\xi+t'\, u) |\, u \rangle- t \langle \nabla F(\xi+t\, u) |\, u \rangle\\
 						& \ge \langle \nabla F(\xi+t'\, u) |\, u \rangle(t-t') +  t'\langle \nabla F(\xi+t'\, u) |\, u \rangle- t \langle \nabla F(\xi+t\, u) |\, u \rangle\\
						& = t  \langle \nabla F(\xi+t'\, u)-  \nabla F(\xi+t\, u)|\, u \rangle\\
						& =  \frac{t}{t'-t}  \langle \nabla F(\xi+t'\, u)-  \nabla F(\xi+t\, u)|\, (\xi+t'\, u-(\xi+t\, u) \rangle>0,
 \end{align*}
 \normalsize
 where we used in the second line the convexity of $F$. The same computations show that $\phi_{_F}(\xi, \xi+t'u)- \phi_{_F}(\xi, \xi+tu)<0$ if $t<t'\le 0$.
 
 \noindent $(b)$  Let $\bar x \!\in \hat \partial U$ and $\eta =\frac 12  \big(d(\xi,\bar x)\wedge 1\big)\!\in (0, 1]$ ($d(\xi,\bar x)=|\xi -\bar x|$ if $\bar x\neq \infty$). Let $x_n$, $n\ge 1$, be a sequence of elements of $U$ such that $x_n  \to \bar x \!\in \widehat \partial U$. For large enough $n$, $x_n \!\in B(\xi, \eta)^c$.  Let $x_{n,\eta}$ be such that $|\xi-x_{n,\eta}| = \eta$ and $x_{n,\eta} \!\in (\xi, x_n)= \{t \xi + (1-t)x_n, t\!\in [0,1]\}$ (geometric interval). It exists by continuity of the norm. Then, by $(a)$ applied with $u_n = \frac{x_n-\xi }{| x_n-\xi|}$, it follows that  $\phi_{_F}(\xi, x_n)\ge \phi_{_F}(\xi, x_{n,\eta})$. Consequently, 
 \[
 \liminf_n \phi_{_F}(\xi,x_n)\ge \liminf _n  \phi_{_F}(\xi,x_{n,\eta})\ge \inf_{y:|\xi-y|=\eta}\phi_{_F}(\xi,y)>0
 \]
 since $\phi$ is uniformly continuous on the compact set $\{\xi\}\times S(\xi, \eta)$. As a consequence
 \[
  \liminf_{x\to \bar x, x\in U}\phi_{_F}(\xi,x)\ge  \inf_{y:|\xi-y|=\eta}\phi_{_F}(\xi,y)>0. 
 \]
 
 \noindent $(c)$ is obvious. 
 \hfill$\Box$
 
\medskip
\noindent {\bf Remarks.} $\bullet$ The functions $x\mapsto \phi_{_F}(\xi, x)$ are generally not convex as can be easily seen e.g. by computing their Hessian when $F$ is three time differentiable. 

\medskip
\noindent $\bullet$ The same monotonicity property holds for the  functions $t\mapsto \phi_{_F}(x+t v,x)$ when $v\!\in S(0,1)$ since 
\[
\partial_t \phi_{_F}(x+t v,x)= \langle  \nabla F(x+tv)-\nabla F(x)\,|\, v\rangle = \frac {1}{t} \langle  \nabla F(x+tv)-\nabla F(x)\,|\, x+tv-x\rangle 
\]
which changes its sign when $t$ crosses $0$ owing to the strict convexity of the function $F$. But, in contrast it is clear  by its very definition that $\xi \mapsto \phi_{_F}(\xi,x)$ is convex for  every $x\!\in U$.
 
 \medskip
  \noindent {\bf P7.} {\em  $F$-Bregman Voronoi cells.}
By analogy with the Voronoi partitions of a quantizer $\Gamma\subset \R^d$ associated w.r.t. a norm, we define below their 
counterparts with respect to the $F$-Bregman divergence as follows. This analogy is stronger than expected (see below) due to the following elementary equivalences.

\color{black}For every $\xi$, $a$, $b\!\in U$, $a\neq b$
\begin{align}
 \label{eq:ineqphi1}
  \phi_{_F}(\xi,a) \begin{pmatrix}< \\
\le\end{pmatrix} \phi_{_F}(\xi,b) &\Longleftrightarrow \langle \nabla F(b)-\nabla F(a) | \xi - a\rangle \begin{pmatrix}< \\
\le\end{pmatrix}\phi_{_F}(a,b) \\
\label{eq:ineqphi2} & \Longleftrightarrow \varphi_{a,b}(\xi) :=  \langle \nabla F(b)-\nabla F(a) | \xi -a\rangle \begin{pmatrix}< \\
\le\end{pmatrix} \phi_{_F}(a,b).
\end{align}

Note that, as $F$ is strictly convex, $  \langle \nabla F(b)-\nabla F(a) | b-a \rangle>0$ so that $\nabla F(a)\neq \nabla F(b)$. Hence  the  linear form $\varphi_{a,b}$ is non zero  since $\varphi_{a,b}(b)>0$. As a consequence, all the sets $\{\varphi_{a,b}= c\}$, $c\!\in \R$, are hyperplanes that partition $\R^d$ into two half-spaces $\{\varphi_{a,b}\le c\}$ and  $\{\varphi_{a,b}>c\}$.

\begin{dfn}[Bregman--Voronoi partition of a quantizer]\label{def:BregVorGrid} Let $\Gamma\subset  U$ be a  finite set. A Borel partition $\big(C_a(\Gamma)\big)_{a\in \Gamma}$ of $U$ satisfying 
 \[
C_a(\Gamma) \subset  \bigcap_{b\in \Gamma} \big \{\xi\!\in U: \phi_{_F}(\xi,a) \le  \phi_{_F}(\xi,b)\big\}
\]
is called a Bregman--Voronoi partition of $U$. 
\end{dfn}
It is clear that for every $a\!\in \Gamma$, $a\!\in C_a(\Gamma)$.  It follows from~\eqref{eq:ineqphi1} and~\eqref{eq:ineqphi2}  that
\[
  \bigcap_{b\in \Gamma, b\neq a } \big \{\xi\!\in U: \phi_{_F}(\xi,a) <  \phi_{_F}(\xi,b)\big\} = \mathring{C}_a(\Gamma)\subset \bar C_a(\Gamma)=  \bigcap_{b\in \Gamma} \big \{\xi\!\in U: \phi_{_F}(\xi,a) \le  \phi_{_F}(\xi,b)\big\}
\] 
(where interior and closure are taken w.r.t. to the  induced topology on $U$) with 
\[
\mathring{C}_a(\Gamma)= \hskip-0.25 cm   \bigcap_{b\neq a,\, b\in \Gamma} \big\{ \xi\!\in U: \varphi_{a,b} (\xi) < \phi_{_F}(a,b)\big\} \; \mbox{ and }\; \bar C_a(\Gamma)=  \bigcap_{b\in \Gamma} \big\{\xi\!\in U: \phi_{a,b} (\xi)    \le \varphi_{_F}(a,b)\big\}.
\]
\color{black}
A crucial fact is that, like for regular Voronoi diagrams  with respect to an Euclidean norm, $\mathring{C}_a(\Gamma)$ and $\bar C_a(\Gamma)$ are  polyhedral  convex sets (open and closed w.r.t. the induced  topology of $U$ respectively). Note that   interior, closure and boundary of the  cells $C_a(\Gamma)$ do not depend on the initial choice of the  Bregman--Voronoi  partition.

\medskip
We associate to any such Bregman--Voronoi partition of $\Gamma$ a  $\phi_{F}$-nearest neighbour projection 
$$\xi \mapsto  {\rm Proj}_{\Gamma}^F(\xi) := \sum_{i=1}^n x_i \textbf {1}_{C_i(x)}\textcolor{black}{(\xi)}
$$ 
and we define the $\phi_{_F}$-quantization  $\widehat X^{\Gamma}$ of a random vector $X:(\Omega, {\cal A})\rightarrow U$ resulting from such a $\phi_{_F}$-Bregman--Voronoi partition by setting as expected 
\[
\widehat X^{\Gamma}= {\rm Proj}_{\Gamma}^F(X)=  \sum_{i=1}^n x_i \textcolor{black}{\textbf {1}_{C_i(x)}(X)}.
\] 

\medskip
We will now adapt for convenience Definition~\ref{def:BregVorGrid}  of the Bregman--Voronoi partition of a grid/quantizer $\Gamma\subset \R^d$ to the  case where we consider an $n$-tuple $x= (x_1,\ldots,x_n)\!\in U^n$. The geometrical way to do that can be simply to say that $\big(C_i(x)\big)_{i=1,\ldots,n}= \big(C_a(\Gamma^x)\big)_{a\in \Gamma^x}$ but this does not allow for a precise definition of $C_i(x)$ and $C_j(x)$ when $x_i= x_j$. 

\begin{dfn}[Bregman--Voronoi partition of an $n$-tuple]\label{def:BregVorx}Let $x\!\in  U^n$. A Borel partition  $\big(C_i(x)\big)_{1\le i \le n}$ of  $U$
satisfying:

\smallskip
 $(i)$ for every $i\!\in \{1,\ldots,n\}$,
\[
C_i(x) \subset  \bigcap_{j=1}^n \big \{\xi\!\in U: \phi_{_F}(\xi,x_i) \le  \phi_{_F}(\xi,x_j)\big\},
\]

 $(ii)$ if $C_i(x)\neq \varnothing$ and $x_j=x_i$, then $C_j(x)= \varnothing$,

\smallskip
 \noindent is called an {\em $F$-Bregman--Voronoi partition} (or more simply a {\em Bregman--Voronoi partition}) induced by the $n$-tuple $x$.  
 \end{dfn}
 
 The $\phi_{F}$-nearest neighbour projection $\xi \mapsto \sum_{i=1}^n x_i \textbf {1}_{C_i(x)}$ resulting from such a $F$-Bregman--Voronoi partition  of $x$ coincides with that    resulting from $\Gamma^x$, in particular
 \begin{equation}\label{eq:quantizerx}
 \widehat X^x:= \sum_{i=1}^n x_i \textbf{1}_{C_i(x)}(X)= \sum_{a\in \Gamma^x}  \textcolor{black}{a} \textbf{1}_{C_a(\Gamma^x)}(X) = \widehat{X}^{\Gamma^x}.
 \end{equation}


The choice to enhance the $n$-tuple rather than the grids in this section is due to the method of proof that we will adopt for the existence of $(r, \phi_{_F})$-optimal quantizer which strongly relies on the  $(r,P,F)$-Bregman distortion functions and analytical properties.

\medskip
 Now, it follows again from~\eqref{eq:ineqphi1}-\eqref{eq:ineqphi2}   that any Bregman--Voronoi  partition of an $n$-tuple $x\!\in U^n$  with  {\em pairwise distinct components} satisfies
$$
\mathring{C}_i(x)= \hskip-0.25cm \bigcap_{\,x_j\neq x_i} \big \{\xi\!\in U: \phi_{_F}(\xi,x_i)< \phi_{_F}(\xi,x_j)\big\} \; \mbox{ and }\;  \bar{C}_i(x) =  \bigcap_{j=1}^n \big \{\xi\!\in U : \phi_{_F}(\xi,x_i) \le  \phi_{_F}(\xi,x_j)\big\},
$$
where interior and closure should still be understood with respect to (the induced  topology on) $U$.
$\bar C_i(x)$  
are  {\em polyhedral convex subsets of $U$} that cover $U$ and whose intersections are contained in finitely many traces of hyperplanes on $U$.

\bigskip
\noindent {\bf Remark.} Again, if   $F$ and $\nabla F$ can be continuously extended to $\bar U$ (closure in $\R^d$), then the above partitions can be extended to the closure  $\bar U$ of $U$ as well as the interior and the closure of the cells (relatively to $\bar U$). The (finite) quantizer $\Gamma$  (resp. the $n$-tuple ) can also be taken in $\bar U$ (resp. in $\bar U^n$).

\section{Optimal  quantization with respect to Bregman divergence  $\phi_{_F}$ ($r=2$)}\label{sec:optiBregquant}
$\rhd$ {\em The optimal  quantization problem for the Bregman divergence}. Solving this problem for the Bregman divergence $\phi_{_F}$  induced by $F$ first consists in proving that the function $G_{r,n}$ attains its {\em minimum} over $U^n$. Such an $n$-tuple,  solution to
\[
 {\rm argmin}_{U^n}G_{r,n},
\]
 if any, is called an $(r, \phi_{_F})$-optimal quantizer of $P$ (or $X$) at level $n$ for  the Bregman divergence $F$ or, in short, an $(r, \phi_{_F})$-optimal $n$-quantizer of $P$ (or $X$). 
 



 \subsection{Optimal $\phi_{_F}$-quantization  at level  $n=1$ ($r=2$)} 
 
The following proposition  (mainly Claim~$(b)$) turns out to be   a slight generalization   of~\cite [Proposition~2.1]{Fischer2010} when $n=1$. 
\begin{prop}\label{prop:existence} 
Assume the distribution $P$ of $X$ satisfies the  moment assumption~\eqref{eq:Integrabilite}.
%
%

\smallskip
\noindent  $(a)$ The distortion function  $G_{1}:U\to \R_+$ attains a unique strict minimum on $U$ at $\E\, X\!\in U$.

\smallskip
\noindent $(b)$ 
%
Furthermore the (l.s.c.)  extension   $G_{1}:  \bar U^{\widehat{\R^d}}\to \R_+$ still attains a unique strict minimum on $U$ at $\E\, X\!\in U$. When $d=1$ the same holds for the l.s.c. extension  of $G_1$ to $ \bar U^{\bar \R}$.
\end{prop}

\noindent {\em Proof.} $(a)$ Note that  
\[
G_{1}(x)= \int_{U}\phi_{_F}(\xi,x) P(\mathrm{d}\xi) = \E\, F(X)- F(x)- \langle \nabla F(x)\,|\,\E\, X-x\rangle  
\]
can be rewritten
\begin{equation}\label{eq:G_1}
G_{1}(x) = \E \, F(X) -F(\E\,X) + \phi_{_F}(\E\,X,x)
\end{equation}
so that 
\begin{equation}\label{eq:minn=1}
\min_U G_{1} = \E\, F(X)-F(\E\,X) \;\mbox{ and } \;   \displaystyle {\rm argmin}_{U}\, G_{1} =\{\E\,X\} \subset U
\end{equation}
since $\phi\geq 0$,  \textcolor{black}{$\phi_{_F}(u,v)=0$ if and only if $ u=v$ by the strict convexity of $F$}.

\smallskip 
\noindent $(b)$  
%
 Let $\bar x \!\in\bar U^{\widehat{\R^d}}$. One has by l.s.c. of $G_1$ on $U^{\widehat{\R^d}}$ that 
 \begin{align*}
 G_1(\bar x) & \le \liminf_{x\to \bar x, x\in U}\E\, \phi_{_F}(X,  x)\\
 		   & =\E\, F(X)-F(\E\,X) +  \liminf_{x\to \bar x, x\in U}\E\, \phi_{_F}(\E\, X, x)\\
		   & =   \E\, F(X)-F(\E\,X) +\phi_{_F}(\E\, X,\bar x)
 \end{align*}
by definition~\eqref{eq:lscextphi}  of $\phi_{_F}(\xi,\cdot)$ on the Alexandroff boundary of $U$. But, as $\E\, X \!\in U$, it follows from Lemma~\ref{lem:monotonie}$(b)$ that  $\phi_{_F}(\E\, X,\bar x)>0$. Consequently, noting  $\min_{U}G_1= G_1(\E\, X) =  \E\, F(X)-F(\E\,X)$ completes the proof. 
\hfill$\Box$

\medskip
\noindent {\bf Remark.} $\bullet$ Note claim~$(a)$ can be extended {\em mutatis mutandis} to the case where $U$ is possibly not open but simply convex.

\smallskip
\noindent $\bullet$  The quantity 
\[
\min_UG_1 = \E\, F(X)-F(\E\, X)
\]
can be seen as a kind of $\phi{_F}$-variance or an $F$-variance of $X$.

\subsection{Optimal $\phi_{_F}$-quantization  at levels $n\ge 2$ ($r=2$)}\label{subsec:existr=2}
%
 One can decompose the distortion  $G_{r,n}(x)$ on any Borel $F$-Bregman--Voronoi partition as follows
\[
G_{n}(x_1,\ldots,x_n)  = \sum_{i=1}^n\int _{C_i(x)}\phi_{_F}(\xi,x_i) P(\mathrm{d}\xi).
\]
When $G_n$ can be extended on $\bar U^n$~--~e.g. because $F$ and $\nabla F$~can themselves be extended on $\bar U$~--~one easily checks that the above decomposition still holds true.
\paragraph{Existence and properties of an optimal $F$-Bregman $n$-quantizer.}

 The existence result below is a slight generalization of former  existence results (see~\cite{Banerjeeetal2005} and an improvement and extension~\cite[Proposition~2.1]{Fischer2010} (this  paper being mainly focused on Bregman divergence based quantization in  Hilbert and Banach spaces). The proof below follows the strategy by induction developed in~\cite{Pages1998,PagCEMRACS} or~\cite{LuPag23} for   $L^r$-Voronoi quantization associated to a norm on $\R^d$ or on a Hilbert space (see also~\cite{GrafL2000} or~\cite{CuestaM1998} for other approaches). 


\begin{thm}[Existence for $r=2$ and $n\ge 2$] \label{thm:Existence}   Assume that the distribution $P$ satisfies the moment assumption~\eqref{eq:Integrabilite} (with $r=2$). Moreover, if     
 $U$ is unbounded and $d\ge 2$, assume that  the $F$-Bregman divergence $\phi= \phi_{_F}$ satisfies
\begin{equation}\label{eq:boundary-lsc}
\forall \xi \!\in U
, \quad \liminf_{x\to \infty,\, x\in U}\phi_{_F}(\xi, x)= \sup _{x\in U} \phi_{_F}(\xi,x).
\end{equation} 


\medskip
\noindent $(a)$ Then for every $n\ge1$ there exists an $n$-tuple $x^{(n)}= \big(x^{(n)}_1, \ldots, x^{(n)}_n\big) \!\in U^n$ which minimizes $G_{n}$ over $U^n$. Moreover, if  the support (in $U$) of the distribution $P$ has at least $n$ points then $x^{(n)}$ has pairwise distinct components and $P\big(\mathring C_i(x^{(n)})\big)>0$ for every $i\!\in \{1,\ldots,n\}$.

\smallskip
\noindent $(b)$ The distribution  $P$ assigns no mass to the boundary of Bregman--Voronoi partitions of $x^{(n)}$, i.e.
\[
P\left (\bigcup_{i=1}^n\partial C_i\big( x^{(n)})\right) = 0,
\]
and $x^{(n)}$ satisfies the $\phi_{_F}$-stationary (or master) equation
\begin{equation}\label{eq:Mastereq}
P\big( C_i(x^{(n)})\big)\,x^{(n)}_i  -  \int_{ C_i(x^{(n)})} \xi P(\mathrm{d}\xi) =0, \quad  i =1,\ldots,n.
\end{equation}
 
\noindent $(c)$ The sequence $\displaystyle G_{n}\big( x^{(n)}\big) = \min_{U^n} G_{n}$ decreases as long as it is not $0$ and converges to $0$ as $n$ goes to $+\infty$.

\noindent $(d)$ When $d=1$, the above claims $(a)$-$(b)$-$(c)$ remain true without assuming~\eqref{eq:boundary-lsc}.
\end{thm}

\noindent {\bf Practitioner's corner} ({\em Shrinking $U$ may help}). A priori, $U$ has been defined as the domain of  definition of the function $F$ but, it is clear that if  $P$ is in fact supported by a smaller open convex set $V$, one can apply the above theorem to $V$, the restriction $F_{|V}$ of $F$ to $V$ and the  resulting restriction of the Bregman divergence $(\phi_{_F})_{|V\times V}= \Phi_{|F_{V}}$ to $V\times V$. With, as a by-product that the optimal $n$-quantizers in $x^{(n)}$ will be $V^n$-valued since then $V$ plays the  role of $U$ in the above theorem. Thus, Condition~\eqref{eq:boundary-lsc}  may be easier to satisfy e.g. if $V$ is bounded while $U$ is not. This fact can also be applied when $r>2$ in the next Section.

\bigskip
\noindent {\bf Remark}. 
$\bullet$ The assumption made on $F$ seems slightly lighter than those made in the literature (see e.g.~\cite{Fischer2010},~\cite{Chatterji1973}) since we only ask the functions $\phi_{_F}(\xi,\cdot)$ to ``attain'' their supremum at $\infty$ (if $U$ is unbounded) and not at any point of $\partial U$. 

\noindent $\bullet$   Actually, Condition~\eqref{eq:boundary-lsc} at $\infty$   also reads \textcolor{black}{for every $\xi \!\in U$}, $\displaystyle  \liminf_{x\to \infty, \,x\in U} \phi_{_F}(\xi, x) =   \sup_{x\in U} \phi_{_F}(\xi,x)$. It  is in fact equivalent to the seemingly stronger 
\[
\forall\,  \xi\!\in U, \quad  \lim_{x\to \infty, \,x\in U} \phi_{_F}(\xi, x) =    \sup_{x\in U} \phi_{_F}(\xi,x).
\]
It  is of course satisfied if, for every $\xi\!\in U$, $\displaystyle \lim_{x\to \infty, \,x\in U}  \phi_{_F}(\xi, x)=+\infty$. 
Such is the case if 
\begin{equation}\label{eq:critere}
\lim_{x\to \infty, \,x\in U}  \langle \nabla F(x)\,|\, x\rangle -F(x) = +\infty.
\end{equation}

\begin{cor} If $x\!\in U^n$ is a stationary quantizer at level $n\ge 1$  (i.e. has pairwise distinct components with $P$-negligible Bregman--Voronoi boundary and  satisfies the $\phi_{_F}$-master equation~\eqref{eq:Mastereq}), then all induced quantizations $\widehat X$ are $\P$-$a.s.$ equal and
\begin{equation}
e_{2}(x, P,\phi_{_F})^2 = \E \, F(X)-\E\, F\big(\widehat X^{x}\big).
\end{equation}
In particular, if $x^{(n)}$ is an optimal $F$-Bregman quantizer then
 \begin{equation}\label{eq:magique}
e_{2,n}(P,\phi_{_F})^2 = \E \, F(X)-\E\, F\big(\widehat X^{x^{(n)}}\big).
\end{equation}
\end{cor}

\noindent {\em Proof.} If $x\! \in U^n$  is stationary  then~\eqref{eq:Mastereq} reads
\[
\forall\, i\!\in \{1, \ldots,n\},\quad \E\big(X\,|\,X\!\in C_i(x) \big) =x_i 
\]
\textcolor{black}{with the convention that if $P\big(C_i(x)\big)=0$ then we set  $ \E\big(X\,|\,X\!\in C_i(x) \big) =x_i$. Hence}, by Formula~\eqref{eq:quantizerx} for $\widehat X^x$ based on  the nearest  $\phi_{_F}$-neighbour projection
\[
\E\, (X\,|\, \widehat  X^x) =\widehat X^x.
\]
As a consequence
\begin{align*}
e_2(x, P,\phi_{_F})^2 &= \E\, \phi_{_F}(X,\widehat X^x)\\
				&= \E\big(F(X)-F\big( \widehat X^x\big) -\langle \nabla F(\widehat X^x)\,|\, X-\widehat X^x\rangle \big) \\
				& = \E \Big(F(X)-F\big(\widehat X^x\big)  -\E \big(\langle \nabla F(\widehat X^x)\,|\, X-\widehat X^ x\rangle \,|\, \widehat X^x \big)\Big) \\
				&=  \E\big(F(X)-F\big( \widehat X^x\big)\big) -\E\Big(\E\big(\langle \nabla F(\widehat X^x)\,|\, X-\widehat X^x\rangle \,|\, \widehat X^x\big)\Big)\\
				& =  \E\big(F(X)-F\big( \widehat X^x\big)\big) -\E\Big(\langle \nabla F(\widehat X^x)\,|\, \underbrace{ \E\big( X-\widehat X^x\rangle \,|\, \widehat X^x\big)}_{=0}\Big)
\end{align*}
Finally, one gets this simple formula
\begin{equation}\label{eq:magique-Intro}
		e_2(x, P,\phi_{_F})^2		 = \E\big(F(X)-F\big( \widehat X^x\big)\big).
\end{equation}
The application to $\phi_{_F}$-optimal  quantizer is straightforward.
\hfill$\Box$

\bigskip
\noindent {\bf Remark.} The above result can be rewritten with grids $\Gamma$ as follows
\[
\E\big(X\,|\,\widehat X^{\Gamma}\big)= \widehat X^{\Gamma}.
\]

%
%
%

\noindent {\bf Practitioner's corner.} By analogy with optimal Voronoi  quantization theory, Equation~\eqref{eq:Mastereq} can  be called the {\em master equation} of Bregman--Voronoi optimal $\phi_{_F}$-quantization since it is the starting point of all methods to  compute optimal quantizers. In the above form~\eqref{eq:Mastereq}, it appears in gradient descents algorithms (a.k.a.  Competitive Learning Vector Quantization  for regular quantization).  When implementing Lloyd's fixed point procedure (also known as $k$-means in unsupervised learning and clustering), one rewrites it in its original form  (see 
\cite{PagesY2016}), i.e. as a fixed point equation (where $X\sim P$): 
\begin{equation}\label{eq:Mastereq'}
x^{(n)}_i  =  \frac{\int_{ C_i(x^{(n)})} \xi P(\mathrm{d}\xi)}{P( C_i(x^{(n)}))} = \E\,\big(X\,|\, X\in C_i(x^{(n)})\big)\!\in  C_i(x^{(n)}), \quad =1,\ldots,n.
\end{equation}

\medskip
\noindent {\em Proof of Theorem~\ref{thm:Existence}}. $(a)$ The proof is based on an induction on the quantization level $n$. We assume that both divergences $\phi_{_F}(\xi, \cdot)$ and distortion   $G_n$ are extended to $\bar U^{\widehat{\R^d}}$ and $\big( \bar U^{\widehat{\R^d}}\big)^n$ by~\eqref{eq:lscextphi} and ~\eqref{eq:lscGrn} (their l.s.c. envelopes) respectively. \textcolor{black}{Our aim is to prove that these extensions attain a minimum in $U^n$ i.e. that ${\rm argmin}_{(\bar U^{\widehat{\R^d}})^n}G_n \subset U^n$ and is nonempty}.

\smallskip
\noindent $\rhd$  $n=1$. This follows from  
Proposition~\ref{prop:existence}$(b)$.

\noindent$\rhd$  $n\Rightarrow n+1$. Let $\displaystyle x^{(n)}\!\in {\rm argmin}_{U^n}\, G_n=  {\rm argmin}_{\bar U^n}\, G_n\neq \varnothing$ supposed to be  nonempty by the induction assumption. 

If $G_n(x^{(n)})=0$, then $|{\rm supp}_U(P)|\le n$ and $\min_{U^{n+1}}G_{n+1} = G_{n+1}(x^{(n)},x^{(n)}_1)= G_n(x^{(n)})=0$.

If $G_n(x^{(n)})>0$, then \textcolor{black}{$|{\rm supp}_U(P)|\ge n+1$ since $x^{(n)}$ has pairwise distinct components}. Moreover,  there exists $\displaystyle x_{n+1}\!\in {\rm supp}_U(P)$ such that $P\big(B(x_{n+1},\eta)\cap U  \big) >0$ for every $\eta>0$. 
Then,  $\displaystyle \big\{\xi\!\in U: \phi_{_F}(\xi,x_{n+1})< \min_{i=1,\ldots,n} \phi_{_F}(\xi,x^{(n)}_i)\big\}$ is  an open set  since   $\displaystyle \xi\mapsto \phi_{_F}(\xi, x_{n+1})-\min_{i=1,\ldots,n}  \phi_{_F}(\xi, x^{(n)}_i)$ is  continuous on $U$,  nonempty   since  it contains $x_{n+1}$. Consequently it has a positive 
$P$-probability which clearly entails that $G_{n+1}\big(x^{(n)},x_{n+1}\big)< G_n(x^{(n)}) = \min_{U^n} G_n$. 

%

Set $\lambda_{n+1} =  G_{n+1}(x^{(n)},x_{n+1})<  \min_{U^n} G_n  $ and  $B_{n+1} = \big\{x\!\in (\bar U^{\widehat{\R^d}})^{n+1}: G_{n+1}(x)\le \lambda_{n+1}\big\}$. 
The set $B_{n+1}$ is clearly closed in $(\bar U^{\widehat{\R^d}})^{n+1}$ since (the l.s.c.  extension to $\widehat{\R^d}$ of) $G_{n+1}$ is l.s.c. Hence $B_{n+1}$ is compact in $\bar U^{\widehat{\R^d}}$. Consequently, $G_{n+1}$ attains its minimum on $B_{n+1}$ at some $n+1$-tuple $x^{(n+1)} = \big(x^{(n+1)}_1, \ldots, x^{(n+1)}_{n+1}\big)$ which is clearly a minimum of $G_{n+1}$ on the whole $(\bar U^{\widehat{\R^d}})^{n+1}$.   Moreover, $x^{(n+1)}$ has pairwise distinct components,  otherwise $G_{n+1}\big( x^{(n+1)}\big) = G_n\big((x^{(n+1)}_i)_{i=1,\ldots,n+1, i\neq i_0} \big)$ for some index $i_0$, hence  $G_{n+1}\big( x^{(n+1)}\big)\ge \min_{U^n} G_n>\lambda_{n+1}$ which is impossible. 

Let us prove that $x^{(n+1)}$ lies in $U^{n+1}$ so  that $x^{(n+1)}\!\in {\rm argmin}_{U^{n+1}} G_{n+1}$. 
First, if ($U$ is unbounded and) $x^{(n+1)}_{i}= \infty$ for some $i\!\in \{1, \ldots,n+1\}$, then,  by Condition~\eqref{eq:boundary-lsc} for every $\xi \!\in U$, $\phi_{_F}(\xi, x^{(n+1)}_{i})\ge \phi_{_F}(\xi,u)$ for every $\! u\in U$. In particular $\phi_{_F}(\xi, x^{(n+1)}_{i})\ge\min_{\ell\neq i}\phi_{_F}(\xi, x^{(n+1)}_{\ell})$  so that $G_{n+1}(x^{(n+1)})= G_n\big( (x^{(n+1)}_\ell)_{\ell\neq i}\big) > \lambda_{n+1}$ which contradicts that $x^{(n+1)}\!\in B_{n+1}$.  Hence $x^{(n+1)}\!\in   \bar U^{n+1}$ ($\bar U$ closure in $\R^d$). 

Assume now there exists an $i$  such that $x^{(n+1)}_{i} \!\in\partial U=  \bar U \setminus U$.  Consider an $\phi_{_F}$-Bregman Voronoi partition $\big( C_i(x^{(n+1)})\big)_{i=1,\ldots,n}$ of $U$ such that $U_i:=C_i(x^{(n+1)})$ is convex and open (in $U$ hence in $\R^d$). This is always possible since $x^{(n+1)}_i\neq \infty$ and $x^{(n+1)}$ has pairwise distinct components. 
Note that  $P(U_i)>0$ since, otherwise, $\int_{ C_{i}(x^{(n+1)})} \phi_{_F}(\xi,x^{(n+1)}_{i} )P(\mathrm{d}\xi)=0$ which would imply that  $G_{n+1}(x^{(n+1)})\ge  G_n\big( (x^{(n+1)}_\ell)_{\ell\neq i}\big) >\lambda_{n+1}$. 

Set $P_i := P(\,\cdot\,|\, U_i)$. Then it follows from Proposition~\ref{prop:existence}$(b)$ that
the  function $g_i(v):= \int_{U_i} \phi_{_F}(\xi,v)P(\mathrm{d}\xi)$, $v\!\in  \bar U_i^{\widehat{\R^d}}$ attains its minimum on $\bar U_i^{\widehat{\R^d}}$ only at $\widetilde x_i = \int_{U_i} \xi P_i(\mathrm{d}\xi)\!\in U_i$.  Hence $g_i(x^{(n+1)}_i) >g_i(\widetilde x_i)$ so that 
%
\begin{align*}
G_{n+1}\big( x^{(n+1)}\big) &= \sum_{j\neq i}  \int_{U\setminus   C_j(x^{(n+1)})}  \phi_{_F}(\xi,x^{(n+1)}_j)P(\mathrm{d}\xi) + \int_{ C_i(x^{(n+1)})}  \phi_{_F}(\xi,x^{(n+1)}_i)P(\mathrm{d}\xi) \\
				   &> \sum_{j\neq i}  \int_{U\setminus   C_j(x^{(n+1)})}  \phi_{_F}(\xi,x^{(n+1)}_j)P(\mathrm{d}\xi) + P(U_i)\int_{C_i(x^{(n+1)})}  \phi_{_F}(\xi,\tilde x_i)P_i(\mathrm{d}\xi) \\
				   &= \sum_{j\neq i}  \int_{U\setminus   C_j(x^{(n+1)})}  \phi_{_F}(\xi,x^{(n+1)}_j)P(\mathrm{d}\xi) + \int_{C_i(x^{(n+1)})}  \phi_{_F}(\xi,\tilde x_i)P(\mathrm{d}\xi)\\
				  &\ge G_{n+1}\big( \tilde x^{(n+1)} \big),
\end{align*}
where $\tilde x^{(n+1)}_j=  x^{(n+1)}_j$, $j\neq i$ and $\tilde x^{(n+1)}_i= \tilde x_i$. This violates the minimality of $x^{(n+1)}$ so that $g_i(x^{(n+1)}_i) \le g_i(\widetilde x_i)$ which implies in turn that $x^{(n+1)}_i= \widetilde x_i \!\in U_i= C_{i}(x^{(n+1)})\subset U$ and  
\[
x^{(n+1)}_i = \int_{U_i}\xi P_i(\mathrm{d}\xi) = \frac{\int_{ C_{i}(x^{(n+1)})} \xi P(\mathrm{d}\xi)}{P( C_{i}(x^{(n+1)}))}.
\]

\noindent $(b)$ Assume  $P\big(\mathring C_i(x^{(n)})\big) = 0$ for some $i\!\in \{1,\ldots,n\}$.     
\begin{align*}
G_n\big( x^{(n)}\big) &= \int_{U\setminus \mathring C_i(x^{(n)})} \min_{j=1,\ldots,n} \phi_{_F}(\xi,x_j)P(\mathrm{d}\xi)= \int_{U\setminus \mathring C_i(x^{(n)})}\min_{j\neq i} \phi_{_F}(\xi,x_j)P(\mathrm{d}\xi) \\
&=  \int_{U}\min_{j\neq i} \phi_{_F}(\xi,x_j)P(\mathrm{d}\xi)\\
& = G_{n-1}\big(x^{(n)}_j, \, j\neq i\big)\ge \min_{U^{n-1}} G_{n-1}> \min_{U^n} G_{n} 
\end{align*}
 which contradicts the minimality of $x^{(n)}$. Hence $P\big(\mathring C_i(x^{(n)})\big) > 0$. 
 
   Set $U_i = \mathring C_i(x^{(n)})$ and $P_i= P(\cdot \,|\, U_i)$. Then applying~\eqref{eq:minn=1} from Proposition~\ref{prop:existence} (i.e. the case $n=1$), we derive that  $u\mapsto \int_{U_i} \phi_{_F}(\xi,u)P_i(\mathrm{d}\xi)$ has a strict minimum at $x^*_i= \int \xi P_i(\mathrm{d}\xi)$. If $x^*_i\neq x^{(n)}_i$ then
\[
\int_{\mathring C_i(x^{(n)})} \phi\big(\xi, x^{(n)}_i \big)P(\mathrm{d}\xi) > \int_{\mathring C_i(x^{(n)})} \phi_{_F}(\xi, x^*_i) P(\mathrm{d}\xi) 
\]\
which in turn implies 
\begin{align*}
G_n\big(x^{(n)}_1,\ldots, x^{(n)}_{i-1},x^*_i,x^{(n)}_{i+1}, \ldots,x^{(n)}_n\big) &\le G_n(x^{(n)}) \\
&\quad\; -\int_{ \mathring C_i(x^{(n)})} \phi_{_F}(\xi,x^{(n)}_i)P(\mathrm{d}\xi) + \int_{ \mathring C_i(x^{(n)})} \phi_{_F}(\xi,x^{*}_i)P(\mathrm{d}\xi)\\
& < G_n(x^{(n)}).
\end{align*}
This  contradicts the minimality of $x^{(n)}$. Hence $x^{(n)} _i= x^*_i= \frac{\int_{\mathring C_i(x^{(n)})} \xi P(\mathrm{d}\xi)}{P( \mathring C_i(x^{(n)}))}$ or, equivalently,
\begin{equation}\label{eq:master1}
P\big( \mathring C_i(x^{(n)})\big) x^{(n)} _i = \int_{\mathring C_i(x^{(n)})} \xi P(\mathrm{d}\xi), \quad i=1,\ldots,n.
\end{equation}

To alleviate notation let us temporarily drop $x^{(n)}$ in the notation of the Bregman--Voronoi partition under consideration.  

Assume now that there exists $i_0$, $j_0\!\in \{1, \ldots,n\}$ such that $P(\bar C_{i_0}\cap \bar C_{j_0})>0$. We consider a new Bregman--Voronoi partition $(C'_i)_{i=1,\ldots,n}$ satisfying
\[
C'_{i_0} = \mathring C_{i_0}\cup (\bar C_{i_0}\cap \bar C_{j_0}), \; C'_{j_0} = C_{j_0}\setminus (\bar C_{i_0}\cap \bar C_{j_0})
\]
and $C'_i\supset C_i$, $i\neq i_0$, $j_0$. Such a diagram exists by assigning in a proper (measurable) way those of the  boundary edges of $C_{i_0}$ not   adjacent to $C_{j_0}$ to its  neighbouring cells.

By the same reasoning as above, we show by the optimality of $x^{(n)}$, after noting that   $C'_{i_0}$ is convex,   that 
\begin{equation}\label{eq:master1bis}
P(C'_{i_0}) x^{(n)} _{i_0} = \int_{ C'_{i_0}} \xi P(\mathrm{d}\xi).
\end{equation}

Subtracting~\eqref{eq:master1} (with $i=i_0$) to~\eqref{eq:master1bis} yields
\[
x^{(n)} _{i_0} =\frac{ \int_{ \bar C_{i_0}\cap \bar C_{j_0}} \xi P(\mathrm{d}\xi)}{P(\bar C_{i_0}\cap \bar C_{j_0}) }\!\in \bar C_{i_0}\cap \bar C_{j_0}
\]
since  $\bar C_{i_0}\cap \bar C_{j_0}$ is convex. For every $\xi\!\in  \bar C_{i_0}\cap \bar C_{j_0}$,  $\phi_{_F}(\xi, x^{(n)} _{i_0})= \phi_{_F}(\xi, x^{(n)} _{j_0})$. Setting $\xi = x^{(n)} _{i_0}\!\in \bar C_{i_0}\cap \bar  C_{j_0}$ yields 
\[
\phi_{_F}(x^{(n)} _{i_0},x^{(n)} _{j_0} ) = \phi_{_F}(x^{(n)} _{i_0},x^{(n)} _{i_0} ) = 0
\]
which implies $x^{(n)} _{i_0} =x^{(n)} _{j_0} $.  Hence  a contradiction to the optimality of $x^{(n)}$ which would have at most $n-1$ pairwise distinct components. This proves  that,  for every $i\!\in \{1,\dots,n\}$, $P\big( \partial C_i\cap C_j\big) =0$ for every $j\neq i$ so that $P\big( \partial C_{i}\big) =0$ since $\bar C_i$ is a polyhedral convex set.

\smallskip
\noindent $(c)$ Let $(z_n)_{n\ge 1}$ be an $U$-valued sequence, everywhere dense in $U$. Note that $\phi_{_F}(u, v)\to 0$ as $v\to u$ for every $u\!\in U$ so that, for every $u\!\in U$, $\min_{1\le i\le n}\phi_{_F}(u,z_i)\downarrow 0$ as $n\to +\infty$ since there exists a subsequence $n(u)$ such that $z_{n(u)}\to u$. As  $\displaystyle \int_U \min_{1\le i\le n} \phi_{_F}(u,z_1) P(du) <+\infty$ owing to~\eqref{eq:Integrabilite},  $\displaystyle G_n(z_1, \ldots,z_n) = \int_U \min_{1\le i\le n} \phi_{_F}(u,z_i) P(du) \downarrow 0$ by the Beppo Levi monotone convergence theorem.

\smallskip
\noindent $(d)$ When $d=1$, we can consider the compactification of $\R$ by $\bar \R= [-\infty, +\infty]$ and revisit the proof of $(a)$ using the specific l.s.c. envelope of $\phi_{_F}(\xi,\cdot)$ at $+\infty$ and $-\infty$. But in the proof of $(a)$,  in the  case if $x^{(n+1)}_i = + \infty$ or$-\infty$  for some $i\!\in \{1, \ldots, n+1\}$, the cell $C_i(x^{(n+1)})$  is always convex since it is an interval of  $\R$ and has a positive $P$-weight. Consequently, applying the argument used  right after for cells induced by ``finite'' $x^{(n+1)}_i$ can be applied: $v\mapsto \int_{C_i(x^{(n+1)}}\phi_{_F}(\xi,v)P(\mathrm{d}\xi)$ attains its minimum on $C_i(x^{(n+1)})$ (equal to the conditional mean $\widetilde x_i$ of $P$ given $C_i(x^{(n+1)})$ which lies in $C_i(x^{(n+1)})\subset U$ but should also be equal to $x^{(n+1)}_i$ to preserve the status of global  minimizer which yields the same contradiction since it implies  $\pm \infty= x^{(n+1)}_i=\widetilde x_i\!\in U_i$ to be finite. Hence $x^{(n+1)} \!\in \bar U^{n+1}$. The rest of the   proofs  of claims $(a)$-$(b)$-$(c)$ is unchanged and we made no use of the assumption.\hfill$\Box$

  
\medskip
\noindent {\bf Examples of Bregman divergences as similarity measures in one  dimension.} 
The  examples below are all  $1$-dimensional so that Theorem~\ref{thm:Existence} applies without further  conditions on the functions $F$ (beyond strict convexity and differentiability, 
\begin{enumerate} 
\item {\em Regular quadratic similarity measure}. $F(x)= x^2$, $U= \R$,  $\phi_{_F}(\xi,x) = (\xi-x)^2$.

\item {\em Norm--like similarity measure}. $F(x)= x^a$, $a>1$, $U= (0, +\infty)$, and 
$$
\phi_{_F}(\xi,x)= \xi^a +(a-1)x^a- a\, \xi \,x^{a-1}.
$$
\item {\em Itakura--Saito divergence.} $F(x)= -\log (x)$, $U= (0, +\infty)$, and 
$$
\phi_{_F}(\xi,x) = \log \Big(\frac x\xi\Big)+\frac \xi x-1.
$$

\item {\em I-divergence a.k.a. Kullback--Leibler divergence/similarity measure}.  $F(x) = x\log(x)$, $U= (0, +\infty)$ and  
$$
\phi_{_F}(\xi,x)= \xi\Big( \log\Big(\frac \xi x\Big)-1+\frac x\xi \Big).
$$

\item {\em Logistic similarity measure}. $F(x)= x\log x +(1-x)\log(1-x)$,  $U= (0, 1)$, and 
$$
\phi_{_F}(\xi,x)=\xi \log\Big(\frac \xi x \Big)+(1-\xi)\log\Big(\frac{1-\xi}{1-v}\Big).
$$

\item {\em Softplus similarity measure (smooth approximation of the $RELU$ function~(\footnote{For REctified Linear Unit.}))}. $F(x)=F_a(x)= \log(1+e^{ax})/a$, $a>0$, $U= \R$, and 
$$
\phi_{_F}(\xi,x) = \frac 1a\log\Bigg(  \frac{1+e^{a\xi}}{1+ e^{ax}}\Bigg) - \frac{e^{ax}}{e^{ax} +1}(\xi-x).
$$

\item {\em Soft butterfly similarity measure (smooth approximation of the absolute value function)}. $F(x)=F_a(x)= \frac{\log(\cosh(ax))}{a}$, $a>0$, $U= \R$. Then, with obvious notations
\begin{align*}
\phi_{_F}(\xi,x) &= \frac 1a\log\Bigg(  \frac{\cosh(a\xi)}{\cosh(ax)}\Bigg) - \tanh{(ax)}(\xi-x) \\
&= 2\, \phi_{{\rm SoftPlus}_{2a}}(\xi,x).
\end{align*}
The second equality is in fact a consequence of the elementary identity (see Proposition~\ref{prop:bgredivprop1}) 
\[
\frac{\log \cosh(ax)}{a}= 2 \frac{\log(1+e^{2ax})}{2a} -\Big(\frac{\log 2}{a}+x\Big).
\]
\item {\em Exponential loss}. $F_\rho(x)= e^{a x}$,   $a\!\in \R$, $U= \R$, $\phi_{F_\rho}(\xi,x) = \phi_{F_1}(a\, \xi,a\, x)$   where  \linebreak
$$ 
\phi_{F_1}(\xi,x) = e^\xi -e^x-e^v(\xi- x).
$$ 
 \end{enumerate} 

\medskip However, note that  6. (and 7.)  do not fulfill the assumption contained e.g. in~\cite{Fischer2010} to guarantee existence of optimal quantizers for such Bregman divergences, that is 
\begin{equation}\label{eq:Cond2b}
\forall\,  \bar x \!\in \partial \bar U^{\widehat{\R^d}},\; \forall\,  \xi \!\in U, \quad  \liminf_{x\to \bar x, \,x\in U} \phi_{_F}(\xi,   x) =  \sup_{x\in U} \phi_{_F}(\xi,x).
\end{equation}

 \bigskip
\noindent {\bf Examples of Bregman divergences as similarity measures in higher dimension.} 
 \begin{enumerate} 
 \item {\em Regular quadratic similarity measure}. $F(x)= |x|^2$, $x\!\in U= \R^d$.
 \item {\em (Squared) Mahalanobis distance}. $F(x)= x^*Sx$,  $x\!\in U= \R^d$,  \textcolor{black}{$S\!\in {\cal S}^{++}(d,\R) :=$} \textcolor{black}{${\cal S}^{+}(d,\R)\cap GL(d,\R)$} (symmetric and positive definite), and  
 $$
 \phi_{_F}(\xi,x)= (\xi-x)^*S(\xi-x):=|\xi-x|^2_{_S}.
 $$
\item {\em $f$-marginal divergence (additive similarity measure)}. $F(x_1,\ldots,x_d)= \sum_{i=1}^d f(x_i)$, $f$ strictly convex on $V$, is defined on $U= V^d$, and for every $x=(x_1,\ldots,x_d)$, $\xi=(\xi_1, \ldots,\xi_d)\!\in U$,
$$
\phi_{_F}(\xi,x)=  \sum_{i=1}^d \phi_{_f}(\xi_i,x_i).
$$
  \item {\em $f$-marginal divergence (multiplicative similarity measure)}. $F(x_1,\ldots,x_d)= \prod_{i=1}^d f(x_i)$, $f$ strictly convex on $V$, is defined on $U=V^d$, and for every $x=(x_1,\ldots,x_d)$, $\xi=(\xi_1, \ldots,\xi_d)\!\in U$,
  $$
  \phi_{_F}(\xi,x)=  \sum_{i=1}^d \phi_{_f}(\xi_i,x_i)\prod_{j\neq i}\phi_{_f}(\xi_i,x_j).
  $$
  \item {\em Soft max marginal $f$-divergence as similarity measure}.  $F(x_1,\ldots,x_d)=F_{\lambda}(x_1,\ldots,x_d)= \frac{1}{\lambda}\log\Big(\sum_{i=1}^d e^{\lambda f(x_i)}\Big)$, $f$ strictly convex on $V$, is defined on $U= V^d$ and for every $x=(x_1,\ldots,x_d)$, $\xi=(\xi_1, \ldots,\xi_d)\!\in U^d$,
  $$
  \phi_{_F}(\xi,x)= F_{\lambda}(\xi)-F_{\lambda}(x)-  \frac{\sum_{1\le i\le d} f'(x_i)e^{\lambda f(x_i)}(\xi_i-x_i)}{\sum_{1\le i\le d} e^{\lambda f(x_i)}}.
  $$
  \item {\em Soft norm similarity measure}.  $F(x_1,\ldots,x_d)=\sum_{i=1}^d f_\lambda(\xi_i,x_i)$, $f_a(u)= \log\big( \cosh(\lambda u)\big)/\lambda $, $a>0$ (soft butterfly similarity measure), is defined on $U=\R^d$ and for every $x=(x_1,\ldots,x_d)$, $\xi=(\xi_1, \ldots,\xi_d)\!\in U$,
  $$
  \phi_{_F}(\xi,x)= \sum_{i=1}^d \phi_{f_\lambda}(\xi_i,x_i).
  $$

 \end{enumerate}  
 
  \begin{figure}[h!]
    \centering
 \begin{tabular}{ccc}
      \hskip-0,5cm \includegraphics[scale=0.45]{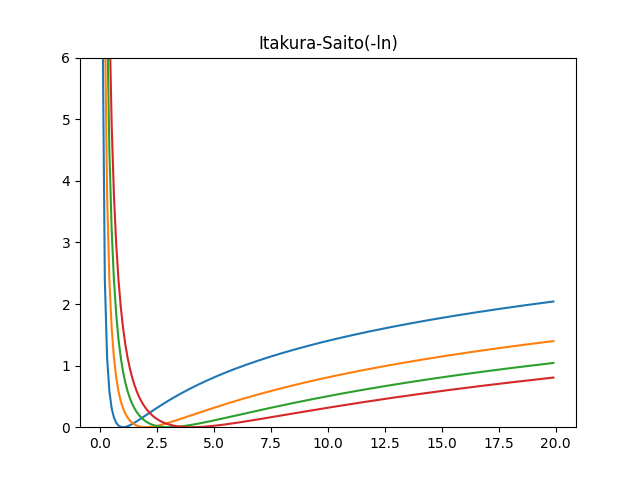}&  \hskip-0,5cm  \includegraphics[scale=0.45]{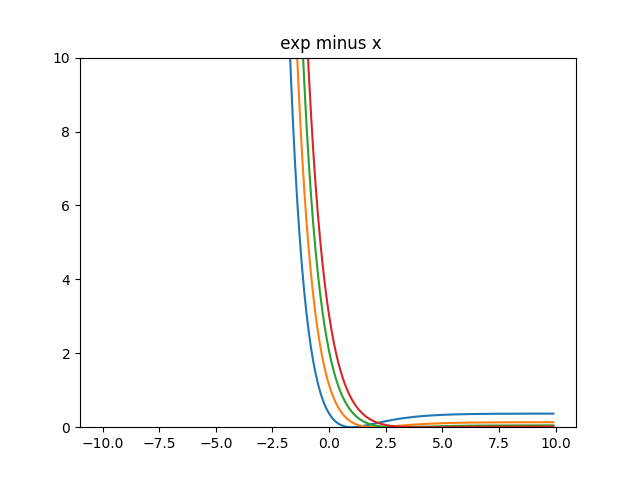}&\\
       \hskip-0,5cm  \includegraphics[scale=0.45]{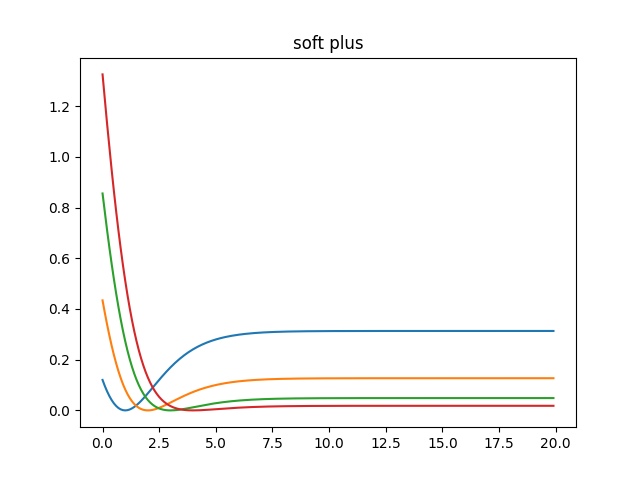}&   \hskip-0,5cm  \includegraphics[scale=0.45]{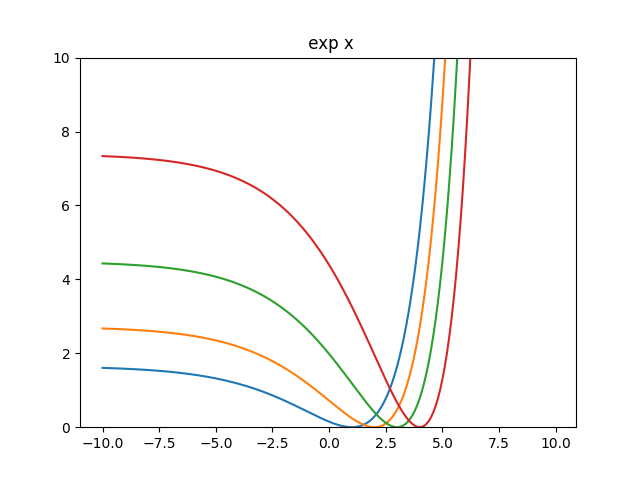}&  
    \end{tabular}
    \caption{\em Non-convexity of $x\mapsto \phi_{_F}(\xi,x)$. $d=1$:  $\xi=1$ (blue), $\xi=2$ (orange), $\xi=3$ (green), $\xi=4$ (purple). Top left: $F(x)=-\log(x)$ (Itakura-Saito divergence); Top right: $F(x)=e^{-x}$ (exponential minus);    Bottom left: $F(x)= \log(1+e^{x})$ (Softplus, $a=1$); Bottom right: $F(x)=e^x$ (exponential).}
    \label{fig:phixinotconvex}
\end{figure}

Figure~\ref{fig:phixinotconvex}  illustrates the   fact  that, in general, the functions $x\mapsto \phi_{_F}(\xi,x)$ $\xi\!\in U$ are not convex. It also illustrates another  fact: for every $\xi\!\in U$, and every vector $u\!\in \R^d$, the $t\mapsto \phi_{_F}(x+tu,x)$ defined in the neighbourhood of   $0$ is decreasing for $t\le 0$ and increasing for $t\ge 0$ since its  derivative is given by the formula
\[
 \partial_t \phi_{_F}(x+tu,x)= -t\langle \nabla F(x+tu,x) -\nabla F(x)|u\rangle
\] 
which has the sign of $t$.

\section{Differentiability of the distortion function  $G^F_{r,n}$ when $F$ is smooth enough}\label{sec:diffdistor}
Another way to retrieve the {\em $\phi_{_F}$-master equation}~\eqref{eq:Mastereq} (when $r=2$) is to compute the gradient of the distortion function $G^F_n$ (see~\eqref{eq:n-distortion}) to write the equation of the critical points of $G^F_n$ since it then  appears as the equation
\[
\nabla G^F_n =0
\]
(provided $\nabla ^2 F$ is positive definite). 

We will see in the next section that, to establish the existence of $(r,\phi_{_F})$-Bregman quantizers,  when $r\neq 2$, this differentiability property will play a key role through the connection between minimality and vanishing of the  gradient.
  
\begin{prop}\label{prop:DiffBregDist} Assume that $F$ is twice differentiable on $U$ and that $P$ satisfies the  moment  assumption~\eqref{eq:Integrabilite}, and assigns no weight to hyperplanes.

\medskip
\noindent $(a)$ {\em Case $r=2$}. Then $G_n$ is differentiable at every  $x\!\in U^n$ having pairwise distinct components and 
\begin{equation}\label{eq:gradientd-dim}
 \nabla G_n (x)= \left[\frac{\partial G_n}{\partial x_i}(x)\right]_{i=1,\ldots,n}={\rm Diag}\big(\nabla^2 F(x_i)\big)_{i=1,\ldots,n} \Big[\int_{C_i(x)}(x_i-\xi)P(\mathrm{d}\xi) \Big]_{i=1,\ldots,n}
\end{equation}
for any Bregman--Voronoi partition $(C_i(x))_{i=1,\ldots,n}$ of $U$ induced by $x$ where ${\rm Diag}\big(\nabla^2 F(x_i)\big)_{i=1,\ldots,n} $ denotes the $dn\times dn$ block diagonal matrix made up with the Hessians of $F$ at the points $x_i$.

\medskip
\noindent $(b)$ {\em Case $r>2$}. Then $G^F_{r,n}$ is differentiable at every  $x\!\in U^n$ having pairwise distinct components and 
\begin{align}
\nonumber \nabla G^F_{r,n} (x)&= \left[\frac{\partial G^F_{r,n}}{\partial x_i}(x)\right]_{i=1,\ldots,n}\\
\label{eq:gradientd-dimr>2} &=\tfrac r2 \big({\rm Diag}\nabla^2 F(x_i)\big)_{i=1,\ldots,n} \left[\int_{C_i(x)}(x_i-\xi)\phi_{_F}(\xi,x_i)^{\frac r2-1}P(\mathrm{d}\xi) \right]_{i=1,\ldots,n}.
 \end{align}
\end{prop}

\noindent{\em Proof.} $(a)$ Let $x\!\in U^n$ with pairwise distinct components. First we note that for every $\xi\!\in U\setminus \cup_{i=1}^n \partial C_i(x)$.
 \[
 \partial_{x_i} \phi_{_F}(\xi,x_i)= \nabla^2F(x_i) (x_i-\xi).
 \]
Let $K\subset U$ be a   convex compact set and let $R_{_K} =\sup_{u\in K} |u|$. Assume  that $x\!\in \mathring{K}^n$. 

Let $u\!\in \mathring{C_i}(x)$. There exists $\eta = \eta_{u,x}>0$ such that for every $\xi\!\in B(x_i,\eta)^n$,
\[
\min_j\phi_{_F}(\xi,x_j) = \phi_{_F}(\xi,x_i). 
\]  
Consequently,  
$$
\partial_{x_i} \min_j \phi_{_F}(\xi,x_j)= \partial_{x_i}\phi_{_F}(\xi,x_i)= -\nabla^2F(x_i)(\xi-x_i).
$$

 If $\xi\! \in \bar C_i(x)^c$, one shows likewise that $\displaystyle \partial_{x_i} \min_j \phi_{_F}(\xi,x_j) =0$.

Moreover, for every $x$, $y\!\in K^n$,
\[
\Big| \min_i \phi_{_F}(\xi,x_i) -\min_i \phi_{_F}(\xi,y_i)   \Big| \le \max_i|\phi_{_F}(\xi,x_i) -  \phi_{_F}(\xi,y_i) |.
\]
Now,
\[
\langle \nabla F(y_i)-\nabla F(x_i)\,|\, \xi-y_i\rangle \le \phi_{_F}(\xi,x_i)-\phi_{_F}(\xi,y_i) \le \langle \nabla F(y_i)-\nabla F(x_i)\,|\, \xi-x_i\rangle
\]
so that 
\begin{align*}
 \max_i|\phi_{_F}(\xi,x_i) -  \phi_{_F}(\xi,y_i) | &\le  \sup_{z \in K^n}\|\nabla^2F(z)\| \textcolor{black}{\max_i (|\xi- x_i| |x_i-y_i|)} \\
 &\le  \sup_{z \in K^n}\|\nabla^2F(z)\|\big(R_{_K} +|\xi|\big) \max_i|x_i-y_i|.
\end{align*}
As $\ds \int_{\R^d}|\xi|P(\mathrm{d}\xi) <+\infty$, one can interchange differentiability and integration which yields the announced result having in mind that the boundaries of the cells $C_i(x)$ are $P$-negligible.

\smallskip
\noindent $(b)$  To keep this technical section short we will not detail the case $r>2$ but rather refer to the proof  of Proposition~\ref{prop:existence_r>2} below which  provides the few specific arguments of this case.\hfill$\Box$ 
 
\section{Optimal $(r,\phi_{_F})$-quantization ($r>2$)}\label{sec:r-opti}
Most results below are established assuming $r\ge 2$ with  different methods of proof    from those developed for the case $r=2$ and, more important, under more stringent assumptions.

 \subsection{Optimal $(r,\phi_{_F})$-quantization at level  $n=1$ ($(r,\phi_{_F})$-Bregman median).}
 

\begin{prop}[$n=1$]\label{prop:existence_r>2} Let $r\!\in (0, +\infty)$.

\smallskip
\noindent $(a)$ Case  $r\ge2$.  Assume that the distribution $P$ of $X$ satisfies the  moment assumption~\eqref{eq:Integrabilite}, that $F$ is $C^2$ on $U$ with $\nabla^2F$ bounded on $U$ and that $\nabla^2F(x)$ is (symmetric) positive definite for  every $x\!\in U$.
%
%
%
%
Assume that at every $\bar x \!\in\partial\,U$, 

\smallskip
\noindent either
\begin{equation} \label{eq:Cond-ra}
 \nabla^2F \mbox{ can  be continuously extended at } \bar x \mbox{ and } \nabla^2F(\bar x) \mbox{ is (symmetric) positive definite}
\end{equation}
or 

\smallskip
 \quad the l.s.c. extensions $\phi_{_F}(\xi, \cdot)$  (on the closure $ \bar U^{\widehat{\R^d}}$ of $U$) satisfy 
\begin{equation}\label{eq:Cond-rb}
 \forall\, \xi \!\in U, \qquad \phi_{_F}(\xi, \bar x) =\displaystyle  \sup_{x\in U}\phi_{_F}(\xi,x)
\end{equation}
and, when $U$ is unbounded, $\bar x = \infty$ satisfies this second condition~\eqref{eq:Cond-rb}.

\smallskip
Then $G_{r,1}:  \bar U^{\widehat{\R^d}}\to \R_+$  attains a unique strict minimum  at  a point $x^{(r)}\!\in U$  \textcolor{black}{which is the unique solution of the stationarity equation
\[
x^{(r)}= \int_U \xi \,\Q_{x^{(r)}}(\mathrm{d}\xi) \!\in U\quad \mbox{ where }\quad \Q_{x} (\mathrm{d}\xi)= \frac{\phi_{_F}(\xi,x)^{\frac r2}\cdot P(\mathrm{d}\xi)}{\int_U\phi_{_F}(\xi,x)^{\frac r2}P(\mathrm{d}\xi)}, \quad x\!\in U.
\]
In particular, $x^{(r)}$ is the  unique  critical point of $G_{r,1}$.}

\noindent \textcolor{black}{$(b)$ Case  $r\!\in (0, 2)$. Assume   all the assumption in $(a)$ are in force and, furthermore, that 
\begin{equation}\label{eq:Integrabilite2}
\exists\, \delta>0,\quad \E\,\phi_{_F}(X,x)^{\frac r2-1-\delta}<+\infty.
\end{equation}
Then,  $G_{r,1} $ attains a minimum on $U$.}
\end{prop}


\smallskip
\noindent  {\bf Remarks.} $\bullet$ Note that when $r=2$ and $n=1$, we proved that the quantizer is  the  expectation $\E\, X$ of $X$ (or the mean of $P$). So, when $n=1$ and $r\neq 2$ the resulting  minimizer of $G_{r,1}$ is a natural candidate to be called an  {\em $(r,\phi_{_F})$-median} of the distribution $P$.

\medskip
\noindent $\bullet$ Note that under the assumptions made on $F$, if $\bar x\!\in \partial U$ satisfies~\eqref{eq:Cond-ra}, then the functions $\phi_{_F}(\xi,\cdot)$ can be extended by continuity at $\bar x$ by the formula
$$
\phi_{_F}(\xi,\bar x)= \int_0^1(1-t) \nabla^2 F(t\xi+(1-t)\bar x) (\xi-\bar x)^{\otimes 2} \mathrm{d}t>0
$$ 
since  $\nabla^2 F$ is positive definite on $U $  and $t\xi+(1-t)\bar x\!\in U$ for every $t\!\in (0,1)$.

 \medskip
\noindent\textcolor{black}{ $\bullet$ The question of uniqueness when $r\!\in (0, 2)$ remains open.}

%

\medskip
\noindent {\em Proof.}  
%
Assume $\bar x\!\in \partial \bar U^{\widehat{\R^d}}$  satisfies~\eqref{eq:Cond-ra} and $\bar x\neq \infty $ (in case $U$  is unbounded) i.e. $\bar x\!\in \partial U$.  As $\nabla^2 F$ is bounded in this neighbourhood of $\bar x$ in $U$, then $\nabla F$ is Lipschitz in the neighbourhood  and can also be continuously extended at $\bar x$.

Let $x\!\in U$.   Then, the function defined on $[0, 1]$ by 
$$
g_x(t) := G_{r,1}\big(tx+(1-t)\bar x\big)
$$ 
is continuous on $[0,1]$. One checks following the lines of Proposition~\ref{prop:DiffBregDist} that as the moment assumption~\eqref{eq:Integrabilite} is in force, $g_x$ is  
differentiable on $(0,1]$  and
\begin{equation}\label{eq:derivgx}
g'_x(t)= \tfrac r2\E \big[ \phi_{_F}(X,tx+(1-t)\bar x)^{\frac  r2-1}\langle \nabla^2F(tx+(1-t)\bar x)(tx+(1-t)\bar x-X), x-\bar x\big\rangle\big].
\end{equation}
Intertwining expectation and differentiation follows from the  fact that, for $y\!\in U\cap \bar B(\bar x,|x-\bar x|)$,
\begin{align*}
|\phi_{_F}(X,y)^{\frac  r2-1}X| & \le \big((1-\tfrac  2r) |\phi_{_F}(X,y)^{\frac  r2}+ \tfrac 2r|X|^{\frac  r2}\big) \\
&\le \tfrac 2r \max\big (\phi_{_F}(X,y)^{\frac  r2}, |X|^{\frac  r2}\big)\le  C_{F,r} \max\big(|F(X)|, |X|)^{\frac  r2}\!\in L^1(\P)
\end{align*}
where we used Young's inequality in the first line. One shows likewise that 
$$
g'_x(t)\to   \frac r2\E \big[ \phi_{_F}(X,\bar x)^{\frac  r2-1}\langle \nabla^2F(\bar x)(\bar x-X), x-\bar x\big\rangle\big] \quad\mbox{as}\quad t\to 0
$$
so that $g_x$ is differentiable at $t=0$ with $g'_x(0)$ formally given by~\eqref{eq:derivgx}. 

%

Moreover, as $X$ is $U$-valued, one has by Taylor's formula with integral remainder that
\[
\phi_{_F}(X, \bar x) = F(X)-F(\bar x)-\langle \nabla F(X)|X-\bar x\rangle = \int _0^1 (1-t)   \nabla^2 F(tX+(1-t)\bar x)(X-\bar x) ^{\otimes 2} \mathrm{d}t >0
\]
since $tX+(1-t)\bar x\!\in U$  and $\nabla^2F(y)$ is positive definite at every $y\!\in U$.  Consequently the  probability measure $P_{\bar x}$ defined by 
$$
P_{\bar x}(\mathrm{d}\xi) = \frac{\phi_{_F}(\xi,\bar x)^{\frac r2-1}}{\E\,[ \phi_{_F}(\xi,\bar x)^{\frac r2-1}]}\cdot P(\mathrm{d}\xi)
$$
is well-defined on $\big(U, {\cal B}or(U)\big)$. Then, we may rewrite $g'_x(0)$ as 
\[
g'_x(0)= - \tfrac r2 \,\E\,[ \phi_{_F}(X,\bar x)^{\frac r2-1}] \cdot \Big \langle \nabla^2 F(\bar x)\Big(\bar x - \int_U \xi P_{\bar x}(\mathrm{d}\xi) \Big)\Big |\bar x- x \Big\rangle.
\] 
 Integrating this equality w.r.t. $P_{\bar x}(dx)$ yields
\begin{align*}
\int g_x'(0)P_{\bar x}(dx) &=  - \tfrac r2 \,\E\,[ \phi_{_F}(X,\bar x)^{\frac r2-1}] \cdot \Big \langle \nabla^2 F(\bar x)\Big(\bar x - \int_U \xi P_{\bar x}(\mathrm{d}\xi) \Big)\Big |\bar x-\int_U \xi P_{\bar x}(\mathrm{d}\xi)  \Big\rangle\\
& =    - \tfrac r2 \,\E\,[ \phi_{_F}(X,\bar x)^{\frac r2-1}] \cdot \nabla^2 F(\bar x)\Big(\bar x - \int_U \xi P_{\bar x}(\mathrm{d}\xi) \Big)^{\otimes 2} \le 0. 
\end{align*}
Assume now that $\bar x$ is a minimizer of $G_{r,1}$ on $\bar U^{\widehat{\R^d}}$. Then $g_x(0) \le g_x(t)$ for  every $t\!\in [0,1]$ so that $g'_x(0) \ge 0$ which in turn implies that $\displaystyle \int g_x'(0)P_{\bar x}(dx) \ge 0$. Hence
$$
 \nabla^2 F(\bar x)\cdot\Big(\bar x - \int_U \xi P_{\bar x}(\mathrm{d}\xi) \Big)^{\otimes 2}  = 0. 
 $$
The Hessian $\nabla^2F(\bar x)$ being positive definite this implies
\[
\bar x = \int_U \xi P_{\bar x}(\mathrm{d}\xi) \!\in U 
\]
which yields a contradiction. 

\smallskip
$\blacktriangleright$ Assume now that $\bar x\!\in \partial \bar U^{\widehat{\R^d}}$  satisfies~\eqref{eq:Cond-rb}. Then, as noted before ,
$$
 \phi_{_F}(\xi, \bar x) = \lim_{x\to \bar x, x\in U}\phi_{_F}(\xi,x)= \sup_{x\in U} \phi_{_F}(\xi,x).
$$ 
Consequently, by Fatou's lemma
 \[
 \E\, \phi_{_F}(X, \bar x)^{\frac r2}= \E\,  \sup_{y\in U} \phi_{_F}(X,y)^{\frac r2}\ge \sup_{y\in U}  \E\, \phi_{_F}(X,y)^{\frac r2} = \sup_{y\in U} G_{r,1}(y).
\]
 Consequently if $\bar x$ is a minimizer then $G_{r,1}$ is constant on $\bar U^{\widehat{\R^d}}$ and, for every $x\!\in U$,  
 \[
  \E\Big(\underbrace{\sup_{x\in U} \phi_{_F}(X,x)^{\frac r2} - \phi_{_F}(X,y)^{\frac r2}}_{\ge 0}\Big)= 0
 \]
 which in turn implies 
 \[
 \forall\, y\!\in U, \quad \phi_{_F}(\xi,y)=\sup_{x\in U}\phi_{_F}(\xi, x)\quad P(\mathrm{d}\xi)\mbox{-} a.s.
 \]
 By continuity of the Bregman divergence $\phi$ on $U\times U$, the functions $\phi_{_F}(\xi,\cdot)$ are all constant hence null since $\phi_{_F}(\xi,\xi)=0$. Which is impossible since $F$ is strictly convex.  
 
 Consequently, any minimizer $x^{(r,1)}$ of $G_{r,1}$ over $\bar U^{\widehat{\R^d}}$ lies in $U$ and satisfies the $(r,\phi_{_F})$-stationarity  or $(r,\phi_{_F})$-master equation
 \[
 x^{(r,1)}= \frac{\int_U \xi\phi_{_F}(\xi,x^{(r)})^{\frac r2-1}P(\mathrm{d}\xi)}{\int_U \phi_{_F}(\xi,x^{(r)})^{\frac r2-1}P(\mathrm{d}\xi)}.
 \]
 
 As for uniqueness, we proceed as follows.  We start from the elementary inequality  that holds for $\rho\ge 1$, 
 \[
 (u+v)^{\rho}\le u^{\rho} +\rho(u+v)^{\rho-1}v, \quad u>0,\; v\ge -u.
 \]
 Then, for every $x,y\!\in U$, it follows, by using~\eqref{eq:phiabc}, that 
 \begin{align*}
 \phi_{_F}(X,x)^{\frac r2}& \le    \phi_{_F}(X,y)^{\frac r2} -\tfrac r2\phi_{_F}(X,x)^{\frac r2-1}  \Big( \phi_{_F}(X,x)-\phi_{_F}(X,y) \Big)\\
 &\le  \phi_{_F}(X,y)^{\frac r2} -\tfrac r2\phi_{_F}(X,x)^{\frac r2-1}  \Big( \phi_{_F}(x,y)-\langle \nabla F (y)-\nabla F(x)\,|\, X-x\rangle \Big),
 \end{align*}
where we used~\eqref{eq:phiabc} in the second line.Taking expectation yields
 \[
 \E\, \phi_{_F}(X,x)^{\frac r2}  \le \E\, \phi_{_F}(X,y)^{\frac r2}  -\tfrac r2\E\, \phi_{_F}(X,x)^{\frac r2-1}\Big( \phi_{_F}(x,y) - \big\langle \nabla F(y)-\nabla F(x)\,|\, \int_U\xi P_x(\mathrm{d}\xi)-x\big\rangle \Big).
 \]
 
 If $x$ is solution of the stationarity equation (in particular if it is a minimizer of $G_{r,1})$, then, for every $y\!\in U$, $y\ne x$ 
 $$
G_{r,1}(x)=\E\, \phi_{_F}(X,x)^{\frac r2}  \le \E\, \phi_{_F}(X,y)^{\frac r2} -\tfrac r2\phi_{_F}(x,y)\E\,[ \phi_{_F}(X,x)^{\frac r2-1}]   < \E\, \phi_{_F}(X,y)^{\frac r2}= G_{r,1}(y).
$$

Hence $x$ is a minimizer of $G_{r,1}$  and if $x'$ is another minimizer, then  $\phi_{_F}(x,x')  =0$ i.e. $x=x'$.This  proves uniqueness. 

\smallskip
\noindent $(b)$ Assumption~\eqref{eq:Integrabilite2} combined with~\eqref{eq:Integrabilite} ensures the differentiability of the function $g_x$. Then the  existence proof of a minimizer of $G_{r,1}$ is unchanged.
\hfill$\Box$

\subsection{Optimal $(r,\phi_{_F})$  quantizers at levels $n\ge 2$}
\begin{thm} \label{thm:Existence-r}   Let $r\!\in [2, +\infty)$. Assume that  the distribution $P$ of $X$ satisfies the  moment assumption~\eqref{eq:Integrabilite}, that  
$F$ is $C^2$ on $U$ with $\nabla^2F(x)$ (symmetric) positive definite for  every $x\!\in U$.
Assume that at every $\bar x \!\in\hat\partial U =  \partial   U$, 

\noindent either
\begin{equation} \label{eq:Cond-ra-bis}
\nabla^2F \mbox{ can  be continuously extended at } \bar x \mbox{ and } \nabla^2F(\bar x) \mbox{ is (symmetric) positive definite}
\end{equation}
or \\
\indent the l.s.c. extensions $\phi_{_F}(\xi, \cdot)$  on $ \bar U^{\widehat{\R^d}}$ satisfy at $\bar x$ 
\begin{equation}\label{eq:Cond-rb2}
 \forall\, \xi \!\in U, \qquad \phi_{_F}(\xi, \bar x) =\displaystyle  \sup_{x\in U}\phi_{_F}(\xi,x).
\end{equation}
Finally,  if $U$ is unbounded, assume that $\bar x = \infty$ satisfies the above condition~\eqref{eq:Cond-rb2}.

\medskip
\noindent $(a)$ Then for every $n\ge1$ there exists an $n$-tuple $x^{(r,n)}= \big(x^{(r,n)}_1, \ldots, x^{(r,n)}_n\big) \!\in U^n$ which minimizes $G_{r,n}$ over $U^n$. Moreover, if  the support (in $U$) of the distribution $P$ has at least $n$ points then $x^{(r,n)}$ has pairwise distinct components and $P\big(\mathring C_i(x^{(r,n)})\big)>0$ for every $i\!\in \{1,\ldots,n\}$.

\smallskip
\noindent $(b)$ The distribution  $P$ assigns no mass to the   boundary of  any $\phi_{_F}$-Bregman--Voronoi partitions of $x^{(n)}$ i.e.
\[
P\left (\bigcup_{i=1}^n\partial C_i\big( x^{(r,n)})\right) = 0
\]
and $x^{(r,n)}$ satisfies the $(r,\phi_{_F})$-stationary (or $(r,\phi_{_F})$-master) equation
\begin{equation}\label{eq:Mastereq-bis}
 x^{(r,n)}_i = \frac{ \int_{C_i(x^{(r,n)})}\xi \phi_{_F}( \xi,x^{(r,n)}_i )^{\frac r2-1}P(\mathrm{d}\xi)}{\int_{C_i(x^{(r,n)})}\phi_{_F}( \xi,x^{(r,n)}_i )^{\frac r2-1}P(\mathrm{d}\xi) } , \quad  i=1,\ldots,n.
\end{equation}
 
\noindent $(c)$ The sequence $\displaystyle G_{r,n}\big( x^{(r,n)}\big) = \min_{U^n} G_{r,n}$ decreases as long as it is not $0$ and converges to $0$ as $n$ goes to infinity.
\end{thm}
 
\smallskip
 \noindent {\bf Remark.} \textcolor{black}{Condition~\eqref{eq:Cond-ra-bis} } implies by standard arguments that $F$ and $\nabla F$ are Lipschitz, locally in the neighbourhood of $\bar x$ in $\partial U$ and, consequently,  can be extended at $\bar x$.
 
\medskip
\noindent {\bf Remark} (Remark ($1D$- setting). In the one-dimensional case, as could be expected,  the Alexandroff compactification of $\R^d$ can be replaced by the compactification $\bar \R = \R \cup\{\pm \infty\}= [-\infty, +\infty]$. Theorems~\ref{thm:Existence-r}  remains true up to assuming that both $\pm\infty$ satisfy~\eqref{eq:Cond-rb}. 

\medskip
\noindent {\em Proof.} The proof is quite similar to that of Theorem~\ref{thm:Existence} and we leave details to the reader.
\hfill$\Box$

\section{A stop at dimension $d=1$: a land of uniqueness for $\log$-concave distributions} \label{sec:BregTrushkin}

For notational convenience we will denote in this section ${\cal S}_n = \{(u_1,\ldots,u_n) \!\in\R^n: u_1<\cdots<u_k<\cdots <u_n\}$,  ${\cal S}_n(a,b)= {\cal S}_n \cap (a,b)^n$ and its closure $\bar {\cal S}_n([a,b])= \bar {\cal S}_n \cap [a,b]^n$.

When $d=1$, the   definition~\eqref{eq:divBregmandef} of the Bregman divergence $\phi_{_F}$ can be extended to closed intervals of $\R$ (of the form $[a,b] \cap \R$, $-\infty \le a <b\le +\infty$) provided $[a,b]\cap \R \subset U$ since then $F$ and $F'$ are  well-defined on this interval.

Thus, the $F$-Bregman boundary between $u$, $v \!\in (a,b)\cap \R$, $u<v$, then reads 
\[
\phi_{_F}(\xi,u) < \phi_{_F}(\xi,v) \Longleftrightarrow \xi <  \varphi(u,v),
\]
where 
\begin{equation}\label{eq:varphi}
 \varphi(u,v):= u + \frac{\phi_{_F}(u,v)}{F'(v)-F'(u)} = v-\frac{\phi_{_F}(v,u)}{F'(v)-F'(u)}.
\end{equation}
(The strict convexity of $F$ implies that $F'$ is increasing hence injective). The function $\varphi$ reads equivalently
\begin{equation}\label{eq:varphi2}
\varphi(u,v) = \frac{F(u)-F(v)-uF'(u)+vF'(v)}{F'(v)-F'(u)}
\end{equation}
so that $\varphi$  can be formally extended on $( (a,b)\cap \R)^2\setminus \{(u,u),\; u\!\in \R\}$ into  a symmetric function. 

Moreover if $F$ is twice differentiable, one can also extend $\varphi$ by continuity on the whole $( (a,b)\cap \R)^2$ by setting $\varphi(u,u) =u$. 

Let ${\cal S}_n( [a,b]) = \big\{(x_1,\ldots,x_n)\!\in ([a,b]\cap \R)^n: a\le x_1\le \cdots\le x_n \le b\big\}$. Consequently for every $x=(x_1,\ldots, x_n)\!\in {\cal S}_n( [a,b])$, $C_i(x)$, $i=1,\ldots,n$ is an interval with endpoints $x^F_{i-\frac 12}$ and $x^F_{i+\frac 12}$ where 
\begin{equation}\label{eq:1Dboundary}
x^F_{\frac 12} = a, \quad x^F_{i+\frac12} = \varphi(x_i,x_{i+1}),\;
\; i=1,\ldots,n-1,\quad x^F_{n+\frac12} = b
\end{equation}
 so that their union is equal to $[a,b]\cap \R$.
\medskip
\begin{prop} Let $a$, $b\!\in \bar \R$ \textcolor{black}{with $a$ or $b$ finite}. Assume there exists an open interval $I\subset \R$ such that $J=[a,b)$, $(a,b]$ or $[a,b]\subset I$ on which the function $F$ is a strictly convex function and subsequently the Bregman divergence $\phi_{_F}$  is well defined as well. Let $\pi: I \to  J$ defined by $\pi(u) =a$ if $u< a$, $\pi(u)=b$ if $u> b$ and $\pi = I_d$ on $J$ (projection from $I$ onto $J$). Assume ${\cal H}\big({\rm supp}(P)\big)= [a,b]\cap \R$ (convex hull of the support of $P$). Then $\phi_{_F}$ and the distortion functions $G_n$ can  be extended on $I^2$ and $I^n$ respectively and
\[
\forall\, x=(x_1,\ldots,x_n)\!\in I^n, \quad G_n\big(\pi(x_1),\ldots,\pi(x_n) \big)\le G_n(x_1,\ldots,x_n). 
\]
In particular, $\inf_{I^n\cap {\cal S}_n}G_n = \inf_{{\cal S}_n([a,b])}G_n$. 
 \end{prop}
 
 Note that this proposition includes distributions $P$ having possibly atoms at the (finite) endpoints of their support.
 
 \medskip
 \noindent {\em Proof.} Let us deal with the case of $a$. Let $x\!\in (-\infty, a]\cap I$ and $u\!\in J$.  One has $\pi(x)=a$ and 
 \[
 \phi_{_F}(u,x)-\phi_{_F}(u,a) = \phi_{_F}(a,x)  + \big(F'(a)-F'(x)\big)(u-a) \ge 0. 
 \] 
 One shows likewise that if $x\!\in [b, +\infty)\cap I$, $\phi_{_F}(u,x)-\phi_{_F}(u,b)\ge 0$. Hence, for every $x\!\in {\cal S}_n \cap I^n$, 
 \[
\hskip 4.25cm  
G_n(x_1, \ldots,x_n) \ge G_n\big( \pi(x_1),\ldots, \pi(x_n)\big).
\hskip 4.5cm \Box
 \]

\subsection{Main results} We will prove the counterpart  of Trushkin's Theorem (originally established for $r=1$ and $r=2$ in~\cite{Trushkin1982} essentially for Euclidean norms)
 in the case $r=2$ for various Bregman divergences.

We assume throughout  this section that the distribution $P$  is {\em strongly unimodal} in the sense that $P= h\cdot \lambda_1$ where $h$ is a $\log$-concave probability density. Then the support of $P$ and $h$ is the closure  of $\{h>0\}= \{\log h> - \infty\}$. Being naturally convex, it is an interval of the form $[a,b]\cap \R$, $-\infty \le a <b\le +\infty$.  Moreover, one shows that the limits $\ds h(a+)$ and $\ds h(b-)$ exist  in $\R_+$. We may assume w.l.g. that 
$$
h(a)=h(a+)\quad\mbox{and}\quad h(b)= h(b-)
$$ 
(with the convention $h(\pm \infty)= 0$ if necessary). Moreover, one checks that, as $h$ is $\log$-concave, $P$ has finite polynomial moments at any order so that the integrability condition boils down to
\begin{equation}\label{eq:integ1D}
\int_a^b |F(\xi)|P(\mathrm{d}\xi)<+\infty.
\end{equation}
Following Theorem~\ref{thm:Existence}$(c)$ there is no additional condition to guarantee the existence of a $F$-Bregman optimal quantizer. Now we are in position to state the Bregman version of Trushkin's uniqueness theorem. 

\begin{thm}[Trushkin's Theorem for $\phi_{_F}$-Bregman quantization ($r=2$)] \label{thm:Bregman-Trushkin}  Assume $P$ is strongly unimodal on $[a,b]\cap \R$ ($a$ and $b$ defined as above) and satisfies~\eqref{eq:integ1D}. Also assume $F$ is 
such that $F \text{ and }F''$ are  positive on $(a,b)$ and either $\log$-concave or $\log$-convex. Then, for every level $n\ge 1$, there is at most one solution to the master equation~\eqref{eq:Mastereq} lying in ${\cal S}_n(a,b)$.
In particular there is exactly one Bregman optimal quantizer.
\end{thm}



%
%
%
%
%
%

 \smallskip
\noindent {\bf Examples of Bregman functions $F$  with $\log$-concave second derivative.} 
\begin{itemize}
\smallskip
\item $F(x)= x^2$, $(a,b)= \R$,  $\phi_{_F}(\xi,x) = (\xi-x)^2$, $U=\R$ and $\log (F'')(x)=  \log 2$ is constant hence concave. 

\smallskip
\item {\em Norm--like}. $F(x)= x^\lambda$, $\lambda>1$, $U= (0, +\infty)$, $\phi_{_F}(\xi,x)= \xi^\lambda +(\lambda-1)x^\lambda- \lambda \xi \,x^{\lambda-1}$ and 
\[
\log F''(x)= \lambda(\lambda-1)\log (x) \quad\mbox{ is  concave.}
\]
%
\item {\em Soft Plus divergence}. $F(x)=F_{\lambda}(x)=  \log(1+e^{\lambda x})/\lambda$, $\lambda >0$, defined on $U=\R$  and 
\[
\log F''(x)=  -2\log(1+e^{\lambda x})+ \lambda x+\log \lambda\quad \mbox{ is concave (and symmetric)}.
\]
The same holds for the Soft butterfly loss function $F(x)=F_{\lambda/2}(x)= 2\,\log(\cosh(\lambda x/2))/\lambda$, $\lambda>0$, as seen before since they differ from an affine function. 

\smallskip
\item {\em Exponential divergence}. $F_\lambda(x)= e^{\lambda x}$, $(a,b)= \R$, $\lambda\!\in \R$ $\phi_{F_\lambda}(\xi,x) = \phi_{F_1}(\lambda \xi,\lambda x)$ with $ \phi_{F_1}(\xi,x) = e^\xi -e^x-e^v(\xi-x)$ and 
$$
\log F''(x)= \lambda x + 2\log \lambda \mbox{ is clearly  \dots affine.}
$$
 
\item {\em I-divergence or Kullback--Leibler divergence}.  $F(x) = x\log(x)$, $(a,b)= (0, +\infty)$, $\phi_{_F}(\xi,x)= \xi\Big( \log\big(\frac \xi x\big)-1+\frac \xi x \Big)$ and 
$$
\log F''(x)= -\log(x)  \mbox{  is  concave.}
$$
\end{itemize}


\medskip
\noindent {\bf Examples of Bregman functions $F$  with $\log$-convex second derivative.} 

\smallskip
\begin{itemize}
\item {\em Itakura--Saito divergence.} $F(x)= -\log (x)$, $U= (0, +\infty)$, $\phi_{_F}(\xi,x) = \log \big(\frac x\xi\big)+\frac \xi x-1$ and 
$$
 \log F''(x) = -2 \log(x) \quad \mbox{ is   convex.}
 $$

\item {\em Logistic divergence}. $F(x)= x\log x +(1-x)\log(1-x)$,  $U= (0, 1)$, $\phi_{_F}(\xi,x)=\xi \log\frac \xi x +(1-\xi)\log\Big(\frac{1-\xi}{1-v}\Big)$ and 
$$
\log F'' (x) =-\log x-\log(1-x) \quad \mbox{ is   convex.}
$$
\end{itemize}

\subsection{The toolbox} Our proof relies on the avatars  of two classical results, namely a refined version of Gershgorin's Lemma and a   version of the Mountain Pass Lemma (which in fact deserves the status of theorem) for functions defined on a compact convex subset of $\R^d$.  

 \begin{lem}[\`A la Gershgorin lemma]\label{lem:Gershgorin}  Let $A= [a_{ij}]_{1\le i,j\le n}$  be an $d\times d$ symmetric tridiagonal matrix and let $L_i =  \sum_{1\le j\le d}a_{ij}$, $i=1,\ldots,d$. Assume that  
\[
\forall\, i \!\in \{1,\ldots, d\}, \quad a_{ii}\ge 0, \quad a_{ii\pm1}< 0, \quad    L_i \ge 0
\]
(with the convention $a_{10}= a_{d d+1}=-1$) and 
$$
L_1\quad \mbox{ or }\quad L_n>0.
$$ 
Then all the  eigenvalues of $A$ are (strictly) positive.
\end{lem}
 
\noindent {\em Proof.}
%
Let $\lambda$ be an eigenvalue of $A$ and $u\!\in\R^n$ an associated eigenvector. Let $i_0\!\in \{1, \ldots,n\}$ be such that $|u_{i_0}|= \max_{1\le i \le n} |u_i|$. We may assume without loss of generality (w.l.g.) that $u_{i_0}=1$. Then
\[
\lambda = a_{i_0i_0} + a_{i_0i_0-1}u_{i_0-1} + a_{i_0i_0+1}u_{i_0+1} \ge  a_{i_0i_0} +a_{i_0i_0-1}|u_{i_0-1}|+ a_{i_0i_0+1}|u_{i_0+1}|\ge L_{i_0}\ge 0.
\]
Moreover, if $\lambda=0$, then the above inequalities  hold as equalities which implies  $L_{i_0}=0$ and  \textcolor{black}{ $u_{i_0i_0+1}= |u_{i_0i_0+1}|=1$ and $u_{i_0i_0-1}= |u_{i_0i_0-1}|=1$}.  One can make the above reasoning with $i=i_0\pm 1$ and carry on by induction. This shows that all the lines have zero sum which is impossible if $L_1$ or $L_n>0$.
\hfill$\Box$

\begin{thm}[Mountain Pass Lemma]\label{thm:MPL}  $(a)$ {\em General case}  (see e.g.~\cite{Kav1993},~\cite{Struwe1990}). If $\Lambda:\R^d\longrightarrow\R_+$ ($\Lambda$ for ``landscape'' when $d\!=\!2$) is
a continuously differentiable function satisfying $\displaystyle \lim_{|x|\to +\infty}\Lambda(x)$  $= +\infty$ and if two distinct zeros of $\nabla \Lambda$ are strict local
minima then $\nabla \Lambda$ has  a third zero  which can be in no case a strict  local minimum.

\smallskip
\noindent $(b)$  {\em  Compact case} (see \cite{LambertonP1996},~\cite{Cohort1998},~\cite[Chapter~8]{LuPag23}). Let $K\!\subset \R^n$ be   a
nonempty compact convex set with a nonempty interior $O$, $\mathring{K} = O$. If $\Lambda :K\longrightarrow \R$ is $C^1$
on $O$, $\nabla \Lambda$ admits a continuous extension  on $K$ satisfying $\{\nabla \Lambda=0\}\!\subset \,O$ 
and if, for every small enough $\varepsilon>0$, $(Id-\varepsilon \nabla \Lambda)(O)\!\subset O$, then $\nabla \Lambda$ has at least a third critical point in $K$ which is not a strict local minimum for $\Lambda$ on $K$. 
\end{thm}

\subsection{Proof of the main uniqueness result  (Theorem~\ref{thm:Bregman-Trushkin})} The method of proof  follows a strategy developed in~\cite{JouPag2021} in a dual quantization framework and in~\cite{LuPag23} for regular $L^r$-Voronoi quantization.  Calling upon the Mountain Pass Lemma  to conclude to uniqueness comes from~\cite{LambertonP1996}. See also~\cite{BouPag1993} devoted to Kohonen's self-organizing maps in one dimension where ``regular'' quadratic quantization appears as a degenerate setting (neighbourhood of each ``unit'' reduced to itself).

\medskip
\noindent {\sc Step~0} {\em Preliminaries}. Any stationary $n$-tuple  lies in ${\cal S}_n(a,b)$ (up to a permutation). In a one-dimensional framework the expression~\eqref{eq:gradientd-dim} of $\nabla G_n$ reads at any point of ${\cal S}_n(a,b)$
\begin{equation}\label{eq:1Dgrad}
\nabla G_n(x) = \Big[ F''(x_i) \int_{x^F_{i-\frac 12}}^{x^F_{i+\frac 12}}(x_i-\xi)P(\mathrm{d}\xi)  \Big]_{i=1,\ldots n}.
\end{equation}

First note  that, if $F''$ admits a continuous extension on $[a,b]\cap \R$,  $\nabla G_n$ can be continuously  extended to  the closure $\bar {\cal S}_n([a,b])$ in $\R^n$  of $\bar {\cal S}_n(a,b)$ since, by $\log$-concavity, $h$ admits such an extension on $[a,b]\cap \R$ (with our conventions).

Moreover we first  check that $\nabla G_{n}$ has no zero on the boundary $ \partial{ \bar {\cal S}}_n[a,b]$.   Thus, temporarily assume that 
 $\nabla G_{n}(x)=0$ and $x_1=a\!\in \R$. Then $\int_{a}^{x^F_{\frac 32}} (a-\xi)h(\xi)d\xi=0$ which implies $x^F_{\frac 32}=a$ since $h>0$ on $(a,b)$. i.e. $x_2=x_1=a$. By induction this implies $x^F_{n+\frac 12}=a$ which is meaningless since $x^F_{n+\frac 12}=b>a$. Hence $x_1>a$. One shows likewise by considering the $i$th component of $\nabla G_{n}$ that one cannot have $x_{i-1}< x_{i}=x_{i+1}$, etc. Finally any critical point of $G_{n}$  lies  in ${\cal S}_n(a,b)$.

\medskip
\noindent {\sc Step~1} {\em  \textcolor{black}{Proof under the  temporary assumption} $a$, $b\!\in \R$, $h(a) +h(b)>0$ and $h$ differentiable on $(a,b)$}.   Note this implies that $P$ has compact support. As $\log F''$ is \textcolor{black}{convex or concave}, it is right differentiable and so is $F''$ by composition with the exponential. 
Then $G_n $ is differentiable at every $n$-tuple $x\!\in {\cal S}_n[a,b]$ and admits a   {\em right} Hessian $\nabla_r^{2} G_n(x)= [a_{ij}(x)]_{1\le i,j\le n}$ given by 
\begin{align*}
a_{ii}(x)&=   F^{(3)}_r(x_i) \int_{x^F_{i-\frac12}}^{x^F_{i+\frac12}} (x_i-\xi)P(\mathrm{d}\xi) + F''(x_i)  \int_{x^F_{i-\frac12}}^{x^F_{i+\frac12}} P(\mathrm{d}\xi) \\
&\quad  -F''(x_i)\Big[ (x^F_{i+\frac 12}-x_i ) h(x^F_{i+\frac 12})\partial_u\varphi(x_i,x_{i+1})\mbox{\bf 1}_{i\neq n}\\
&\hskip 6cm +(x_i-x^F_{i-\frac 12})h(x^F_{i-\frac 12})\partial_v\varphi(x_{i-1},x_{i})\mbox{\bf 1}_{i\neq 1}  \Big] \\
& \hskip 12cm i=1,\ldots,n,\\
a_{ii+1}(x) & =  - F''(x_i)(x^F_{i+\frac 12}-x_i )h(x^F_{i+\frac 12})\partial_v\varphi(x_i,x_{i+1}),\quad i=1,\ldots,n-1,\\
a_{ii+1}(x) &= -F''(x_i)(x_i-x^F_{i-\frac 12} )h(x^F_{i-\frac 12})\partial_u\varphi(x_{i-1},x_{i}),\quad i=2,\ldots,n,
\end{align*}
where the function $\varphi$ is given by~\eqref{eq:varphi} or~\eqref{eq:varphi2}. Note that
\[
\partial_u\varphi(a,\xi ) = \partial_v \varphi(\xi, b) = 0.
\]

Moreover, it should be noticed that one may formally give a sense to $a_{nn+1}(x)$ and $a_{10}(x)$ by the above formulas (using the conventions made on $h(a)$ and $h(b)$). Consequently, if we set $\displaystyle A_i(x)=  \int_{x^F_{i-\frac12}}^{x^F_{i+\frac12}} (x_i-\xi)P(\mathrm{d}\xi) $, $i=1,\ldots,n$, and $\Delta_f = {\rm Diag}(f(x_i), \, i=1,\ldots,n)$ one rewrites the above formula in the more synthetic way 
\[
\nabla^2G_n(x)  = \Delta_{F^{(3)}}(x) \big[A_i(x)]_{i=1,\ldots,n} + \Delta_{F''}(x) \left[ \frac{\partial A_i(x)}{\partial x_j}\right].
\]

If $x$ is a $\phi_{F}$-stationary point, then $A(x) =0$ so that the sum $\Sigma_{L_i}$  of the i$th$ line $L_i$ reads for $i=1,\ldots,n$, 
\begin{align}
\label{eq:SumLi}\Sigma_{L_i}&= F''(x_i)\times\\
\nonumber &\hskip-0.1cm \left( \int_{x^F_{i-\frac12}}^{x^F_{i+\frac12}} \!h(\xi)d\xi -(x^F_{i+\frac 12}-x_i) h(x^F_{i+\frac 12})\psi(x_i,x_{i+1})\mbox{\bf 1}_{i\neq n} - (x_i-x^F_{i-\frac 12})h(x^F_{i-\frac 12})\psi(x_{i-1},x_{i} )\mbox{\bf 1}_{i\neq 1}\!\!\right)\!\!,
\end{align}
where, for every  $\xi, \xi' \!\in (a,b)\cap \R$ 
\begin{equation}\label{eq:psi}
\psi(\xi,\xi') = \partial _u \varphi(\xi,\xi')+\partial_v\varphi(\xi,\xi')= \partial _u \varphi(\xi,\xi')+\partial_u\varphi(\xi',\xi).
\end{equation}
  One checks at this stage  that, if we define 
  $\widetilde \Sigma_{L_i}(x)$ by the above formula~\eqref{eq:SumLi} for $\Sigma_{L_i}$ but where the two indicator functions have been removed then
  $\widetilde \Sigma_{L_i}(x) =  \Sigma_{L_i}(x)$ for $i=2,\ldots,n-1$ and $  \Sigma_{L_i}(x)>  \widetilde \Sigma_{L_i}(x)$  either for $i=1$ or $i=n$ (or both) since $h(x^F_{\frac 12})= h(a) $,   $h(x^F_{n+\frac 12})= h(b) $ and $h(a) +h(b)>0$ by assumption.
  
  Let us show that $\widetilde \Sigma_{L_i} \ge 0$. To this end we first perform an integration by parts, namely
  \begin{align*}
  \int_{x^F_{i-\frac12}}^{x^F_{i+\frac12}} h(\xi)d\xi &= (x^F_{i+\frac12}-x_i)h(x^F_{i+\frac12}) -(x^F_{i-\frac12}-x_i)h(x^F_{i-\frac12}) -\int_{x^F_{i-\frac12}}^{x^F_{i+\frac12}}(\xi-x_i) h_r'(\xi)d\xi,	 
  \end{align*}
 where $h'_r$ denotes the right derivative of $h$ on $(a,b)$ which exists by $\log$-concavity 
 of $h$ and satisfies $h(u)-h(u_0) = \int_{u_0}^{u} h'_r(v)dv$.  As a consequence
  \begin{align}\label{eq:postIPPa}
  \widetilde \Sigma_{L_i}(x)&= F''(x_i) \bigg (\big(x^F_{i+\frac12}-x_i\big)h\big(x^F_{i+\frac12}\big) \big(1-\psi(x_i,x_{i+1})\big) \\
\nonumber &  \quad + \big(x_i-x^F_{i-\frac12}\big) h\big(x^F_{i-\frac12}\big) \big(1-\psi(x_{i-1},x_{i})\big) \\ 
\label{eq:postIPPb}  & \quad + \int_{x^F_{i-\frac12}}^{x_i}(x_i-\xi) \frac{h_r'}{h}(\xi)P(\mathrm{d}\xi)     -  \int_{x_i}^{x^F_{i+\frac12}}(\xi-x_i) \frac{h'_r}{h}(\xi)P(\mathrm{d}\xi)\bigg).
  \end{align}
  
  Since $h$ is $\log$-concave, $\frac {h_r'}{h}$ is non-increasing on $(a,b)$ so that 
  \[
   \int_{x^F_{i-\frac12}}^{x_i}(x_i-\xi) \frac{h_r'}{h}(\xi)P(\mathrm{d}\xi)     -  \int_{x_i}^{x^F_{i+\frac12}}(\xi-x_i) \frac{h_r'}{h}(\xi)P(\mathrm{d}\xi) \ge \frac{h_r'}{h}(x_i) \int_{x^F_{i-\frac12}}^{x^F_{i+\frac12}}(x_i-\xi) P(\mathrm{d}\xi) = 0
  \]
  still using that $x$ is a $\phi_{_F}$-stationary point. This proves that the first term~\eqref{eq:postIPPa} in the formula for $\widetilde \Sigma_{L_i}(x)$ is non-negative.
  
  \smallskip
  \noindent {\sc Step~3} ({\em Proof of $\psi\le 1$ if $F''$ $\log$-convex/concave}). One starts for the identities coming from~\eqref{eq:varphi}. Hence,  for every $u,v\!\in (a,b)$, $u<v$,
  \begin{equation}\label{eq:partialvarphi}
  \partial_u \varphi(u,v)= \frac{F''(u)\phi_{_F}(u,v)}{(F'(v)-F'(u))^2}\quad \mbox { and }\quad  \partial_v \varphi(u,v)= \frac{F''(v)\phi_{_F}(v,u)}{(F'(u)-F'(v))^2}.
  \end{equation}
  Now set 
   \begin{equation}\label{eq:Phiuv}
 \Phi(u,v) = F''(u) \phi_{_F}(u,v) +F''(v)\phi_{_F}(v,u) - (F'(v)-F'(u))^2.
  \end{equation}
 As $F''$ is    positive $\log$-convex or $\log$-concave on $(a,b)$, one deduces that $F''$ is right differentiable with right derivative $F^{(3)}_r$. Elementary   computations yield that for every $u,v \!\in (a,b)$, $u\le v$, the right derivative     $\partial_{u_r}\Phi$ of $\Phi$ w.r.t. $u$ reads
  \[
    \partial_{u_r}\Phi(u,v) = F^{(3)}_r(u)\phi_{_F}(u,v)-F''(u)\phi_{_{F'}}(u,v)
  \]
  and its cross right derivatives  $  \partial^{2}_{u_rv_r}\Phi$ and  $  \partial^{2}_{v_ru_r}\Phi$ of $\Phi$ coincide and read
  \begin{align*}
  \partial^{2}_{u_rv_r}\Phi(u,v)=\partial^{2}_{v_ru_r}\Phi(u,v)& = \big(F_r^{(3)}(u)F''(v) - F_r^{(3)}(v) F''(u)\big) (v-u)\\
  & =   F''(u)F''(v)\Big(  (\log F'')_r'(u) -(\log F'')_r'(v)\Big) (v-u),
  \end{align*}
where we used  that $F'' >0$ on $U=(a,b)$. 

\smallskip
$\blacktriangleright$ Assume  $\log F''$ concave on $(a,b)$. Then it is right differentiable   with a non-decreasing  right derivative $(\log F'')_r'$ on $(a,b)$ so that  
  \[
  \forall\, u,\, v \!\in (a,b), \quad  u\le v,\quad   \partial^{2}_{v_ru_r}\Phi(u,v)\ge 0.
  \]
  This implies that $ v\mapsto  \partial_{u_r}\Phi(u,v)$ is non-decreasing on $[u,b)$ and, in particular, that
  \[
  \partial_{u_r}\Phi(u,v)\ge   \partial_{u_r}\Phi(u,u)=0.
  \]
  As consequence $ \partial_{u_r}\Phi(u,v)\ge 0$ so that $u\mapsto \Phi(u,v)$ is non-decreasing on $(a,v]$. Hence
  \[
  \Phi(u,v)\le \Phi(v,v)=0.
  \]
This finally implies that   the function $\psi$ defined in~\eqref{eq:psi}  satisfies 
\begin{equation}\label{eq:psi<=1}
\forall\, u,\, v, \; u<v, \quad \psi(u,v)\le 1.
\end{equation}
We deduce that the first term~\eqref{eq:postIPPb} in the expression of $\widetilde \Sigma_i(x)$ is non-negative. Finally this shows that
\[
\widetilde \Sigma_{L_i}(x) \ge 0.
\]

$\blacktriangleright$ Assume  now $\log F''$ convex on $(a,b)$. Then $\log F''$  is right differentiable   with a non-decreasing  right derivative $(\log F'')_r'$ on $(a,b)$ so that  
  \[
  \forall\, u,\, v \!\in (a,b), \quad  u\le v,\quad   \partial^{2}_{u_rv_r}\Phi(u,v)\le 0.
  \]
 This implies that $ u\mapsto  \partial_{v_r}\Phi(u,v)$ is non-increasing on $[a,v]$ and, in particular, that
  \[
  \partial_{v_r}\Phi(u,v)\le   \partial_{v_r}\Phi(u,u)=0.
  \]
  As consequence $v\mapsto \Phi(u,v)$ is non-decreasing on $[v,b)$ so that
  \[
  \Phi(u,v)\le \Phi(u,u)=0
  \]
and one concludes as above that $\widetilde \Sigma_i(x) \ge 0$.

 \smallskip
 \noindent {\sc Step~4} {\em First conclusion: the compact case using the Mountain Pass lemma}. We assume  that the additional assumption in force in Step~3  still holds and that $F''$  can be continuously extended on the closed interval $[a,b]\cap \R$. it follows from what precedes that $\Sigma_{L_i}(x) \ge 0$ for $ i=2,\ldots,n-1$ and $\Sigma_{L_i}(x) >0$ for $i=1,2$. Consequently Gershgorin's Lemma (Lemma~\ref{lem:Gershgorin}) implies that any stationary (or critical) point  of $G_n$ is a strictly local minimum.  Moreover, we will prove now that for  small enough $\varepsilon>0$, $I_n-\varepsilon \nabla G_n$ maps ${\cal S}_n(a,b)$ into itself. 
 
 Let  $x\!\in {\cal S}_n(a,b)$ and $\varepsilon >0$. For every $\xi\!\in \mathring{C}_i(x)$ and $\xi' \!\in \mathring{C}_{i+1}(x)$, one has 
  \[
  x_{i+1}-\varepsilon F''(x_{i+1}) (x_{i+1}-\xi') > x_{i+1} - \varepsilon  \big(x_{i+1}-x^F_{i-\frac 12}\big) = x_{i+1} - \varepsilon \frac{F''(x_{i+1})\phi_{_F}(x_{i+1},x_i)}{F'(x_{i+1})-F'(x_i)}
  \]
 where we used~\eqref{eq:varphi} to establish  the equality. Similarly, 
 \[
 x_i- \varepsilon F''(x_i)(x_i-\xi) < x_i +\varepsilon \frac{F''(x_i)\phi_{_F}(x_i,x_{i+1})}{F'(x_{i+1})-F'(x_i)}.
 \]
 Consequently, using the definition of the function $\psi$ from~\eqref{eq:psi},
 \[
 x_{i+1}-x_i > \varepsilon \frac{F''(x_i)\phi_{_F}(x_i,x_{i+1})+F''(x_{i+1})\phi_{_F}(x_{i+1},x_i)}{F'(x_{i+1})-F'(x_i)}= \varepsilon \psi(x_i,x_{i+1}) \big(F'(x_{i+1})-F'(x_i) \big).
 \]
Since $\psi \le 1$ by~\eqref{eq:psi<=1}, the above inequality is \textcolor{black}{implied by} 
 \[
 x_{i+1}-x_i > \varepsilon  \big(F'(x_{i+1})-F'(x_i) \big)
 \]
which
\textcolor{black}{in turn}  holds as soon as $\varepsilon < \frac{1}{\|F''\|_{\sup}}$ (as $F''$  has a continuous extension on the compact interval $[a,b]$, $ \frac{1}{\|F''\|_{\sup}}>0$). 
\color{black}As a consequence, for every $\xi\!\in \mathring{C}_i(x)$ and $\xi' \!\in \mathring{C}_{i+1}(x)$,  
$$
 x_{i+1}-\varepsilon F''(x_{i+1}) (x_{i+1}-\xi') -\big(x_i- \varepsilon F''(x_i)(x_i-\xi) \big)>0
 $$
Integrating w.r.t. the product measure $\mbox{\bf 1}_{ \mathring{C}_i(x)}(\xi)\cdot P(d\xi)\, \mbox{\bf 1}_{ \mathring{C}_i(x)}(\xi')\cdot P(d\xi') $ yields that $(x-\ve \nabla G_n(x))_{i+1}>(x-\ve  \nabla  G_n(x))_{i}$.
It remains to prove that $(x-\ve  \nabla  G_n(x))_1 >a$ and $(x-\ve  \nabla  G_n(x))_n <b$. Note that, for every $\xi\!\in \mathring{C}_1(x)$
\[
x_1-\ve F''(x_1)(x_1-\xi)= x_1 (1-\ve F''(x_1))+\ve   F''(x_1)\xi  
\]
if $\varepsilon < \frac{1}{\|F''\|_{\sup}}$ then $\ve F''(x_1)\!\in (0,1)$ so that the right hand side of the above equality is a convex combination of two real numbers greater than $a$, hence $x_1-\ve F''(x_1)(x_1-\xi)>a$. Integrating w.r.t. $\mbox{\bf 1}_{ \mathring{C}_i(x)}(\xi)\cdot P(d\xi)$ proves that $(G_n)_1(x) >a$. One proceeds likewise for the other inequality. This proves as announced that $I_n-\varepsilon \nabla G_n$ maps ${\cal S}_n(a,b)$ into itself.  
\color{black}

All the assumptions of the Mountain pass Lemma (Theorem~\ref{thm:MPL}) are satisfied, hence  it is  impossible for all critical points of $G_n$  to be strict local minima except if there is only one such critical point. 

\smallskip
\noindent {Step~5} {\em Extension to non-compactly supported  distributions}. Assume there is a level $n\ge 2$ such that $G_n$ has two critical points, say $x^{[1]}$ and $x^{[2]}$, lying in  ${\cal S}_n(a,b)$. Let $K_k= [a_k,b_k]$, $k\ge 1$,  be an increasing sequence of nontrivial compact intervals of $(a,b)$  increasing to $(a,b)$. Let  
$P_k= \frac{h\cdot \mbox{\bf 1}_{K_k}}{P(K_k)}\cdot \lambda_1$, $k\ge 1$, denote the conditional distributions of $P$ given $K_k$. Then, $P_k$ has a $\log$-concave density $h\cdot \frac{\mbox{\bf 1}_{K_k}}{P(K_k)}$, continuous on its support with $h(a_k)$ and $h_k(b_k)>0$. The distortion function $G^{(k)}_{n}$ of $P_k$   can trivially be extended to $\bar {\cal S}_n([a,b])$.

Now, it follows from the definition~\eqref{eq:divBregmandef} of the Bregman divergence $\phi_{_F}$ that 
\begin{align*}
\forall\, \xi\!\in U,\; \forall\, x\!\in K^n,\; \min_{i=1,\ldots,n}\phi_{_F}(\xi,x_i) &\le |F(\xi)| +\sup_{v\in K}|F(v)| + \sup_{y\in K}|F'(v)|\big( |\xi|+\sup_{v\in K} |v|\big)\\
& \le C_{_{K,F}}\big(|F(\xi)]\vee |\xi| +1\big),
\end{align*}
where $C_K= 2 \sup_{v\in K} \big(|F(v)]\vee \big(|F'(v)|(|v|\vee1)\big)\big)$. Then for every $x\!\in K^n$
\begin{align*}
\Bigg|\int_{K_k} \min_{i=1,\ldots,n}\phi_{_F}(\xi,x_i) \frac{P(\mathrm{d}\xi)}{P(K_k)} &-\int_{U} \min_{i=1,\ldots,n}\phi_{_F}(\xi,x_i)  P(\mathrm{d}\xi) \Bigg|\\
&\le \Big(1-\frac{1}{P(K_k)}\Big) \int_{K_k} \min_{i=1,\ldots,n}\phi_{_F}(\xi,x_i)  P(\mathrm{d}\xi) \\
&\hskip 4cm + \int_{U\setminus K_k} \min_{i=1,\ldots,n}\phi_{_F}(\xi,x_i)  P(\mathrm{d}\xi) \\
&\le C_{_{K, P}}\Bigg( \Big(1-\frac{1}{P(K_k)}\Big) \int_{U}\big(|F(\xi)]\vee |\xi| +1\big)P(\mathrm{d}\xi)  \\
&\hskip 3 cm  + \int_{U\setminus K_k}\big(|F(\xi)]\vee |\xi| +1\big)P(\mathrm{d}\xi) \Bigg).
\end{align*}
Hence
\[
\limsup_{k\to +\infty} \sup_{x\in K^n} \Bigg|\int_{K_k} \min_{i=1,\ldots,n}\phi_{_F}(\xi,x_i) \frac{P(\mathrm{d}\xi)}{P(K_k)} -\int_{U} \min_{i=1,\ldots,n}\phi_{_F}(\xi,x_i)  P(\mathrm{d}\xi) \Bigg| =0,
\]
where we used that $P(K_k)\to 1$ as $k\to+\infty$ and $\displaystyle  \int_{U\setminus K_k}\big(|F(\xi)]\vee |\xi| +1\big)P(\mathrm{d}\xi) \to 0$ as $k\to +\infty$ owing to Lebesgue's dominated convergence Theorem. Consequently $G^{(k)}_{n} \to G_n$ uniformly on compact sets of ${\cal S}_n(a,b)$ as $k\to +\infty$.


For large enough $k$, $x^{[1]}$ and $x^{[2]}$ lie in ${\cal S}_n(a_k,b_k)$.    There exists an interval   $B(x^{[1]},s)\subset  {\cal S}_n(a_k,b_k)$ not containing $x^{[2]}$, such that $G_{n}$ is minimum at $x^{[1]}$ by definition of a local minimum. Hence $\min_{B(x^{[1]},s/2)} G_{n}^{(k)}$ cannot be attained on the boundary of $B(x^{[1]},s/2)$ for large enough $k$.  Hence $G_{n}^{(k)}$ has a critical point lying in $B(x^{[1]},s/2)$. Idem for $x^{[2]}$. Hence $G_{n}^{(k)} $ has two critical points which are (strict) local minima which is impossible by Step~3.\hfill$\Box$

%

%

\subsection{More on the ``doubly''-symmetric case ($r=2$)}\label{subsec:symmetry}

In this section we briefly investigate the case where $P$ is symmetric that is $X\stackrel{d}{=}-X$ if $X\sim P$, $X\!\in L^1(\P)$. Numerical experiments show that although the Bregman divergence is not a symmetric function in general,  we will show that, if $F''$ is even, optimal quantizers of such a symmetric distribution are symmetric at every level $n\ge 1$. By symmetric we mean that if $x^{(n)}$ is an optimal quantizer  and is unique  (e.g. because $F''$ and $P=h\cdot \lambda_1$ are symmetric), then 
\[
\forall\, i\!\in \{1, \ldots,n\}, \quad x^{(n)}_{n+1-i} = x^{(n)}_i.
\]
Numerical tests performed with the Bregman-Lloyd algorithm (see~\cite[Chapter 4]{bouto2024}) 
strongly support this claim, e.g. when $X\sim \mathcal{N}(0,1)$ and $F(x) = \log(1+e^x)$ (see Table~\ref{tab:1} with $n= 25$,...). As before we assume that $U=(a,b)$ where $[-b,b]\cap \R$, $b>0$,  is the convex hull of the support of $P$. 

At this stage we are able to prove that fact essentially for $n=1,\,2,\,3$. When $n=1$, we know that $x^{(1)}_1= \E\, X=0$.

For $n=2$, we proceed as follows under the additional assumption that $P(\{0\})=0$. It follows form the $\phi_{_F}$-master equation that a stationary quantizer at level $2$ is of the form
\[
x_1= \frac{\int_{-b}^0 \xi P(\mathrm{d}\xi)}{P((-b,0])}= 2K(0) <0\quad \mbox{ and }\quad x_2= \frac{\int_0^b \xi P(\mathrm{d}\xi)}{P((0,b])}=-2K(0)>0
\]
and that $\varphi(-2K(0)), 2K(0))=0$. 

In fact, one easily checks that if $F''$ is even, then the function $g:(-b,b)\to \R$, 
\[
g(u) :=\varphi(-u,u) = F(-u)+uF'(-u)-F(u)+uF'(u)
\]
satisfies $g'(u) = u(F''(u)-F''(-u))= 0$. Hence $g=g(0)= 0$ i.e. $\varphi(-u,u)=0$, $u\!\in (-b,b)$. As $2K(0)\!\in (-b,0)$, one concludes that 
\[
x^{(2)}= \big(2K(0), -2K(0)\big).
\]
More generally, one has  the following theorem.

\begin{thm} \label{thm:symetrie} $(a)$ Under the  assumptions of Theorem~\ref{thm:Bregman-Trushkin} if,  furthermore, the distribution $P$ is symmetric and $F''$ is even (symmetric), then, for every level $n\ge1$,  the unique  $\phi_{_F}$-optimal quantizer $x^{(n)}$ of $P$  is itself symmetric in the sense that  
\[
\forall\, i \!\in \{1, \ldots,n\}, \quad x^{(n)}_i=  -x^{(n)}_{n-i+1}.
\] 
$(b)$ Assume $P$ is diffuse  and symmetric and $F''$ is symmetric. Then if $x=(x_i)_{i=1,\ldots,n}$ is a stationary quantizer  then   is $\widetilde x= (x_{n_{i+1}})_{i=1,\ldots,n}$ is also a stationary quantizer.
\end{thm}

\noindent {\em Proof.} $(a)$ First note that if $F''$ is even then $F'(u)-F'(0)$ is odd. Let us define for every $u$, $v\!\in (-b,b)$, $u<v$,
\begin{align*}
g(u,v) &= \varphi(u,v)+\varphi(-v,-u)\\
& =  \frac{(F(u)-uF'(u))-(F(v)-vF'(v)))}{(F'(v)-F'(u))^2} + \frac{(F(-v)-vF'(-v))-(F(-u)+uF'(-u))}{(F'(-u)-F'(-v))^2}\\
&= \frac{(F(u)-uF'(u))-(F(-u)+uF'(-u))}{(F'(v)-F'(u))^2} - \frac{(F(v)-vF'(v)))(F(-v)-vF'(-v))}{(F'(v)-F'(u))^2}
\end{align*}
owing to~\eqref{eq:varphi} and the fact that $F'(-u)-F'(-v)= F'(v)-F'(u)$. Now, note that
\[
F(u)-uF'(u)= F(-u)+uF'(-u)
\]
since both sides of the equality  have $-uF''(u)$ as a derivative and are null at $u=0$. Consequently, for every $u$, $v\!\in (-b,b)$, $u<v$,
\[
 -\varphi(u,v)= \varphi(-v,-u).
\]
 
As $x^{(n)}$ is supposed to be optimal, we know that all boundaries of its Bregman--Voronoi diagram are $P$-negligible so that  we can write without ambiguity $\int_c^d g(\xi)P(\mathrm{d}\xi)$ for $\int_I g(\xi)P(\mathrm{d}\xi)$ where $I$ is an interval with endpoints $c$, $d\!\in \{x^{(n),F}_{i+\frac 12}, i=0,\ldots,n+1\}$ (with the usual conventions for $i=0, n+1$).

As $X\stackrel{d}{=}-X$, i.e. $P$ is invariant by $\xi\mapsto-\xi$, one has  for every $i\!\in \{1,\ldots,n\}$
\[
\int_{x^{(n),F}_{i-\frac12}}^{x^{(n),F}_{i+\frac12}}\big(\xi-x^{(n),F}_{i}\big)P(\mathrm{d}\xi)= \int_{-x^{(n),F}_{i+\frac12}}^{-x^{(n),F}_{i-\frac12}}\big(\xi+ x^{(n),F}_{i}\big)P(\mathrm{d}\xi)=  \int_{-x^{(n),F}_{i+\frac12}}^{-x^{(n),F}_{i-\frac12}}\big(\xi-(-x^{(n),F}_{i})\big)P(\mathrm{d}\xi).
\]
It follows from the above identity satisfied by $\varphi$ that 
\[
-x^{(n),F}_{i+\frac12}= -\varphi\big(-x^{(n),F}_{i},-x^{(n),F}_{i+1}\big)=  -\varphi\big(-x^{(n),F}_{i+1},-x^{(n),F}_{i}\big),\quad i=0,\ldots,n
\]
(with the usual convention for $i=0$ and $n$). Setting
\[
\widetilde x^{(n)}_i = x^{(n)}_{n+1-i}, \quad i=1,\ldots,n,
\]
we derive that $\widetilde x^{(n)}$ satisfies the $\phi_{_F}$-master equation~\eqref{eq:Mastereq} for  the distribution $P$  i.e. $\widetilde x^{(n)}$ is an $\phi_{_F}$-stationary quantizer for $P$.

Consequently, under the  assumptions of Theorem~\ref{thm:Bregman-Trushkin}, uniqueness of $(F,P)$-stationary quantizer implies that $\widetilde x^{(n)}= x^{(n)}$.

\smallskip
\noindent $(b)$ The proof is obvious given the computations carried out in $(a)$. 
\hfill$\Box$

\medskip
\noindent {\bf Illustration.} $\bullet$ {\em SoftPlus similarity measure}. Let $F(x)= \log(1+e^{ax})/a$ (with $a>0$) and let $P=\mathcal{N}(0,1)$. One has
\[
\forall\, x\!\in \R \quad F''(x)=a \frac{e^{ax}}{(1+e^{ax})^2}
\]
It is clear that $\log F''(x)= ax-2\log(1+e^{ax})+2\log a$ is concave and one easily checks that $F''(-x)= F''(x)$. Hence both Theorems~\ref{thm:Bregman-Trushkin} and~\ref{thm:symetrie} so that the optimal quantizer is unique and symmetric.
Using the Bregman-Lloyd fixed point procedure (see e.g.~\cite{Banerjeeetal2005}, see also~\cite[Chapter 4]{bouto2024} for a full mathematical   analysis in one dimension),
we can illustrate this result numerically.

\begin{figure}[h!]
    \centering
    \begin{tabular}{c}
   \includegraphics[scale=0.25]{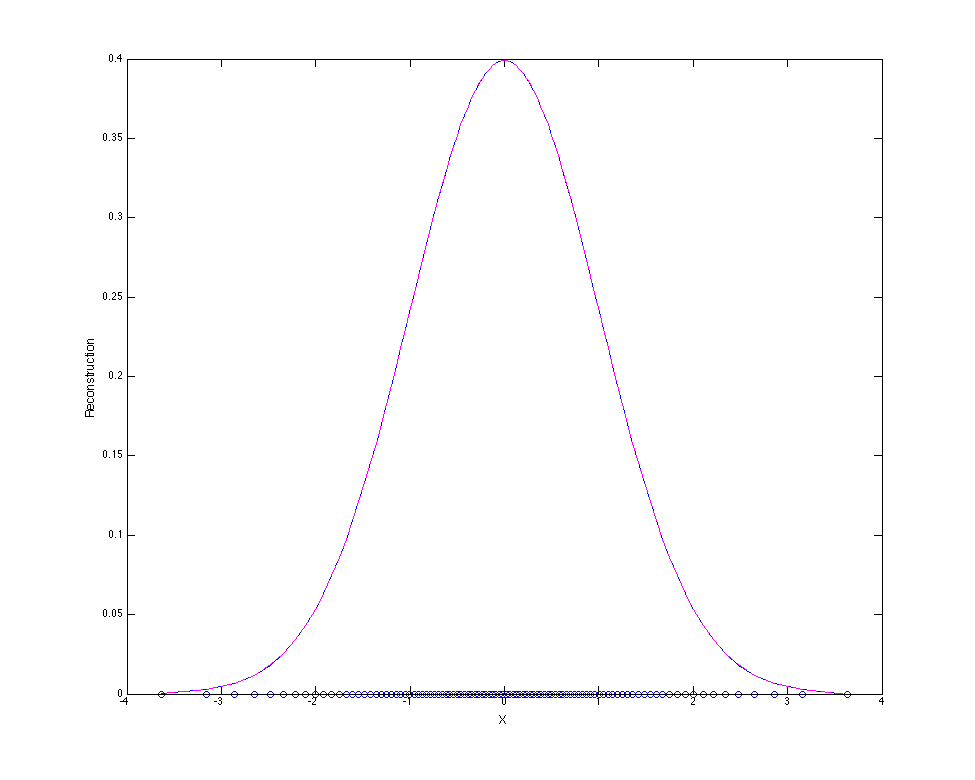}  
    \end{tabular}
    \caption{\em $P(\mathrm{d}\xi) = \frac{e^{-\frac{\xi^2}{2}}}{\sqrt{2\pi}}$. Bregman function $F(x)= \log(1+e^x)$ $(a=1)$.  Level $n = 100$. Quantizer  $x^{(100)}$: \textcolor{blue}{o} on the abscissa axis.  Reconstruction of the density. True density; \textcolor{blue}{---}; reconstructed density \textcolor{magenta}{-- -- --}.}
    \label{fig:1DN010a}
\end{figure}

\begin{figure}[h!]
\centering
\begin{tabular}{c}
   \includegraphics[scale=0.25]{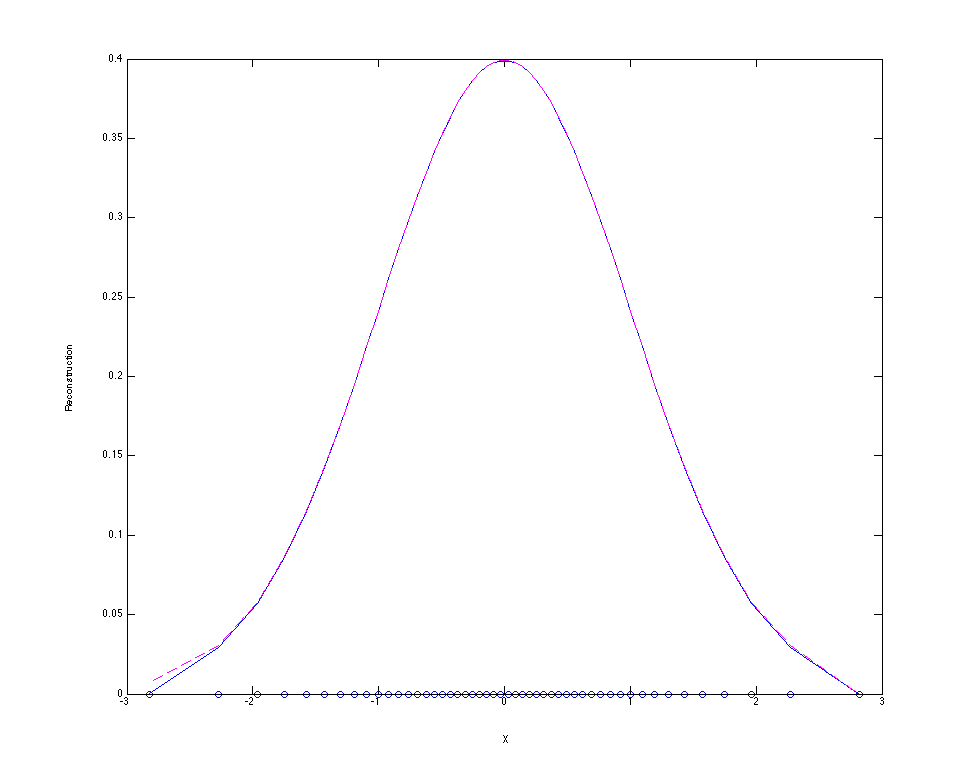}    
   \end{tabular}
    \caption{\em $P(\mathrm{d}\xi) = \frac{e^{-\frac{\xi^2}{2}}}{\sqrt{2\pi}}$. Bregman function $F(x)= \log(1+e^{2x})/2$ ($a=2$).  $n = 100$. Quantizer  $x^{(100)}$: \textcolor{blue}{o} on the abscissa axis. Reconstruction of the density.  True density~\textcolor{blue}{---}; reconstructed density~\textcolor{magenta}{-- -- --}.}
    \label{fig:1DN010b}
\end{figure}

\medskip
See Figures~\ref{fig:1DN010a} and~\ref{fig:1DN010b} of the reconstruction with $n=100$. 

\bigskip

We present in Table~\ref{tab:1} the computed values of an optimal quantizer of  at level (or size)  $n=25$ for the normal distribution $P= \mathcal{N}(0,1)$ with with respect to the Bregman divergences induced by  the two {\em Softplus} functions $F(x)= \log(1+e^{ax})/a$, $a=1$, $2$. 

\begin{table} 
\tiny
\hskip-0,75cm\begin{tabular}{|c||c|c|c|c|c|c|c|c|c|c|c|c|c}
$a$ & & & & & & & & & & & & &\\
\hline
$1$ &  -2.7941  & -2.1936   &-1.8214   &-1.5392  &-1.3053  & -1.1015 &  -0.9177 &  -0.7479  & -0.5882 & -0.4356 
&  -0.2878  &-0.1432  &  0.0000\\
& 0.1432  &  0.2878 &   0.4356  & 0.5882  & 0.7479  & 0.9177 &  1.1015  &  1.3053  & 1.5392  &1.8214   & 2.1936 &   2.7941\\
\hline
$2$ &  -2.4424  & -1.8314  & -1.4897 & -1.2423 &  -1.0434  &-0.8740  &-0.7239 &  -0.5873   &-0.4602 &-0.3399  
&-0.2242  & -0.1114  &  0.0000\\
 &  0.1114 &   0.2242  & 0.3399  & 0.4602   & 0.5873  & 0.7239  & 0.8740  &  1.0434   & 1.2423  &  1.4897  &  1.8314 &   2.4424
\\\end{tabular}

\bigskip
\caption{$n= 25$,  $F(x)= \log(1+e^{ax})/a$, $a=1,2$, $P= N(0,1)$.}
\label{tab:1}
\end{table}
\normalsize

These results may look  at a first glance unexpected since it is clear that the {\em SoftPlus} function itself  is not symmetric. However if one has in mind that the functional $F\longmapsto \phi_{_F}$ is (linear but) not bijective since if $\phi_{_F}$ then $\phi = \phi_{F_{\alpha,\beta}}$ for all functions $F_{\alpha, \beta}(x)= \langle a |x\rangle +\beta$, $\alpha\!\in \R^d$, $\beta \!\in \R$, the paradox is elucidate.  So in this framework  one can replace the  function $F= \hbox{\em SoftPlus}_a$ function by 
\[
F_{1/2, 0}(x) =2 \log(1+e^{x/2})+\tfrac x2
\]
which is symmetric and has the same symmetric second derivative and for which the conclusion of Theorem~\ref{thm:symetrie} would have appeared less misleading but for bad reasons.  Examples of larger optimal $n$-quantizers are proposed in~\cite[Chapter 4]{bouto2024}.

 \smallskip
\noindent $\bullet$  As soon as $F$ is symmetric and ${\cal C}^2$, $F''$ is symmetric as well.  Thus the {\em Soft butterfly} similarity measure  which has a $\log$-concave second derivative satisfies the assumptions of the theorem.


\bigskip
\noindent {\bf Remark.} Theorem~\ref{thm:symetrie} can be extended to other kinds of symmetry properties. If $P$ is invariant by $\xi\mapsto 1-\xi$  and $F''(x)= F''(1-x)$ (e.g. because $F$ satisfies such invariance), then, if $x^{(n)}$ is an optimal quantizer of $P$ at level $n$, then $(1-x^{(n)}_{n+1-i})_{i=1,\ldots, n}$ is also optimal and in case of uniqueness, then $x^{(n)}_i = 1-x^{(n)}_{n+1-i}$, $i=1,\ldots,n$.

\section{Numerical illustrations}

\begin{figure}[h!]
\centering

\vskip-0,25cm \begin{tabular}{cc}
 \includegraphics[scale=0.5]{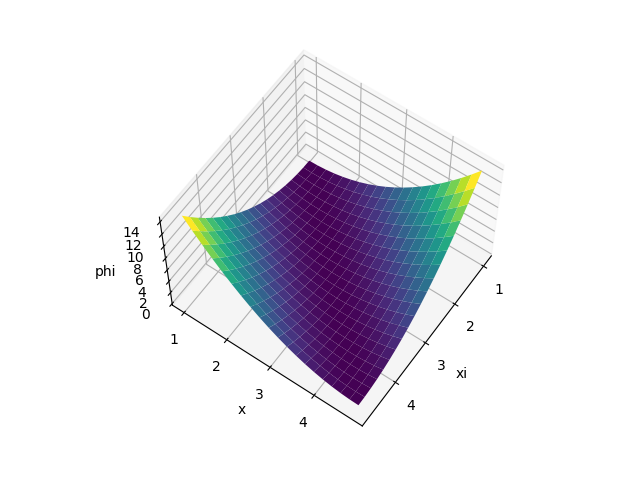}  &   \includegraphics[scale=0.5]{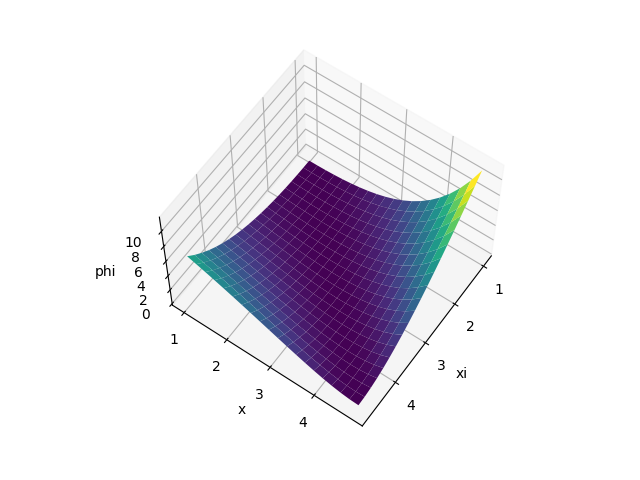}  \\ 
     \includegraphics[scale=0.5]{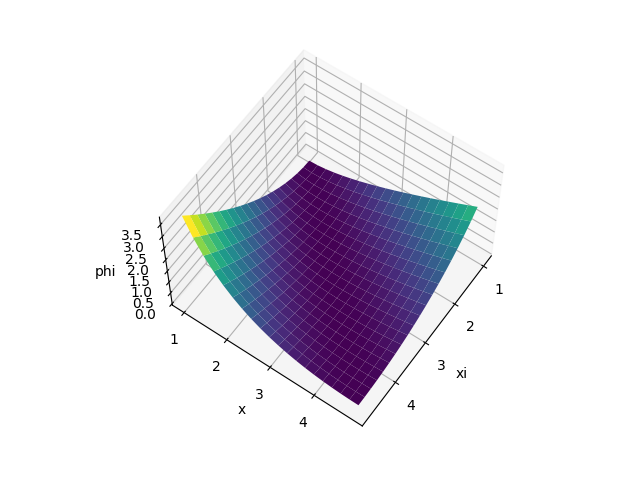} &   \includegraphics[scale=0.5]{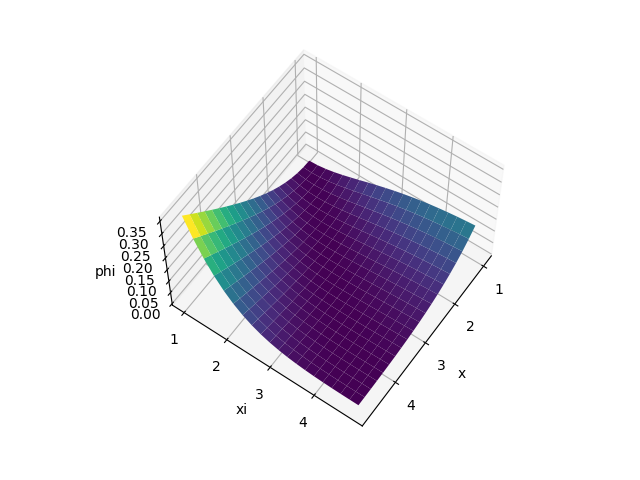}   
  \end{tabular}
    \caption[]{\em $3$-dimensional visualization of Bregman divergences $\phi_F$. From top to bottom and left to right: $F(x)= |x|^2$, $F(x) = e^x$, $F(x)= x\log x$ and $F(x)= \log(1+e^x)$.
\label{fig:Div_3D}}
\end{figure}

$\bullet$ First, we compare in Fig.~\ref{fig:Div_3D} the Bregman divergences variation on the square $[1,4]^2$  for various Bregman divergences $\phi_F$ where $F(x)= |x|^2$ (Squared Euclidean norm), $F(x) = e^x$ (exponential), $F(x)= x\log x$ (Kullback-Leibler) and $F(x)= \log(1+e^x)$ (SoftPlus).

$\bullet$ As a second illustration, we give examples of (quadratc) optimal quantizations of the $2$-dimensional Gaussian distribution  $P= {\cal N}\left(\begin{pmatrix}\tfrac 12 \\ 1 \end{pmatrix};\tfrac{1}{4}I_2\right)$ with a quantization size $n= 41$ w.r.t.  the Bregman divergencs associated to the six functions below. For the last two, , we considered a truncated version of the above Gaussian distribution.

\begin{enumerate}
\item $F(x)=|x|^2$ (regular optimal quadratic quantization).
\item $F(x)=\sum_{i=1,2}\log(1+e^{ax_i})/a$ (SoftPlus, additive marginals, $a=1$)
\item $F(x) =\sum_{i=1,2} e^{a x_i}$ (exponential, additive marginals) $a= 1$ and $-1$.
%
\item $F(x)=\sum_{i=1,2} x_i\log x_i$ (Kullback--Leibler) 
\item $F(x) = -\sum_{i=1,2}\log x_i$ (Itakura--Sa\"ito) 
\end{enumerate}

\medskip
\begin{figure}[h!]
\centering
\begin{tabular}{cc}
  \includegraphics[scale=0.5]{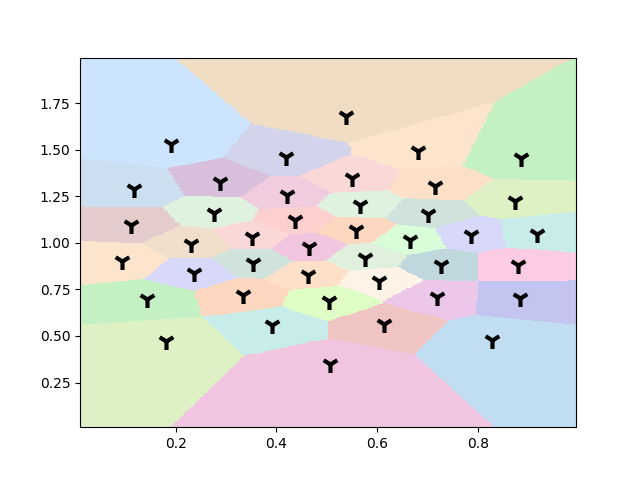}  &   \includegraphics[scale=0.5]{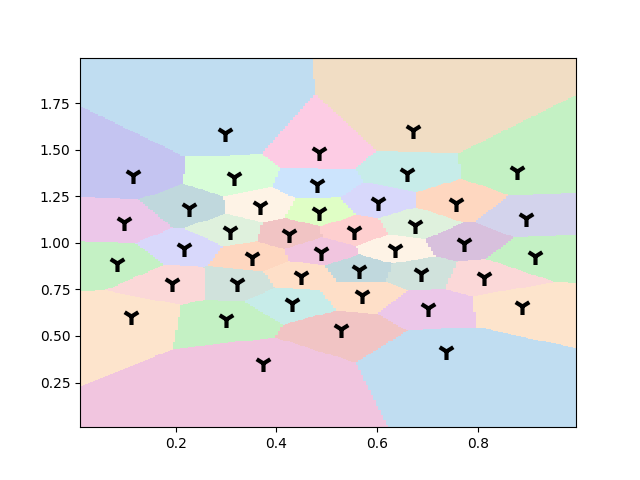}  \\ 
     \includegraphics[scale=0.5]{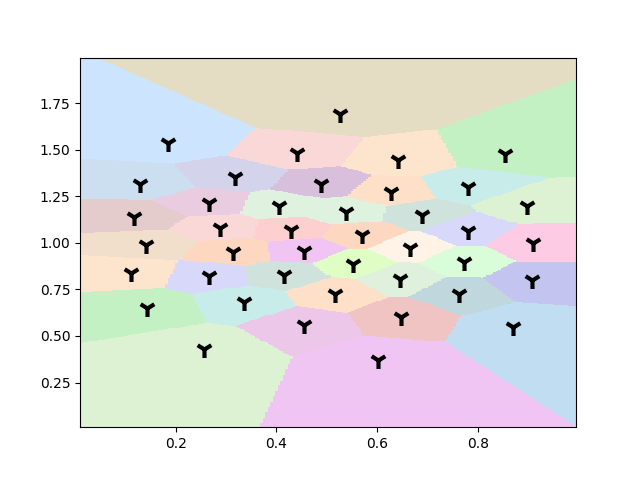} &   \includegraphics[scale=0.5]{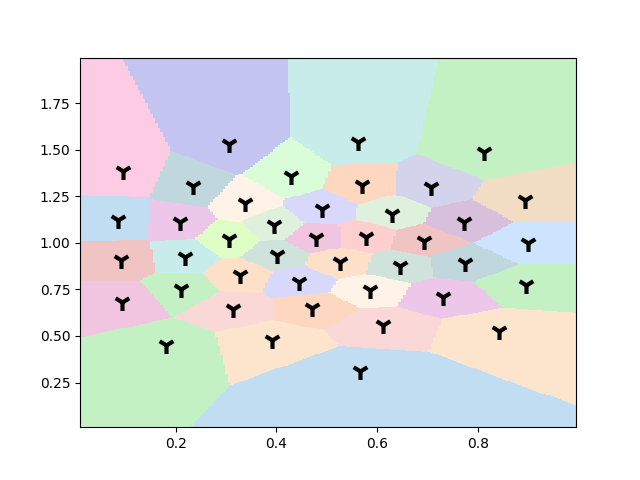} \\
        \includegraphics[scale=0.5]{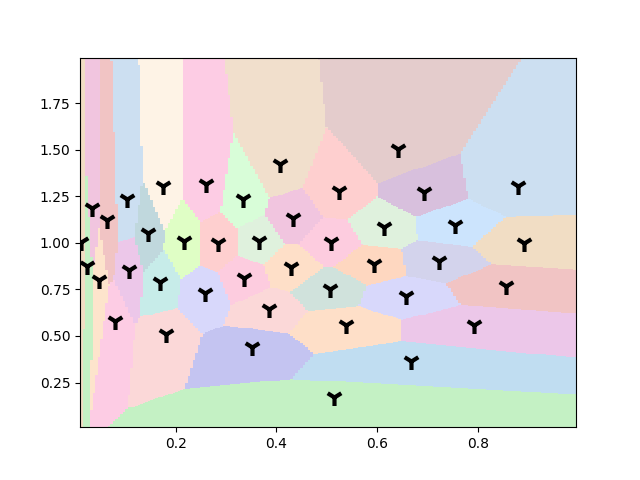}  &  \includegraphics[scale=0.5]{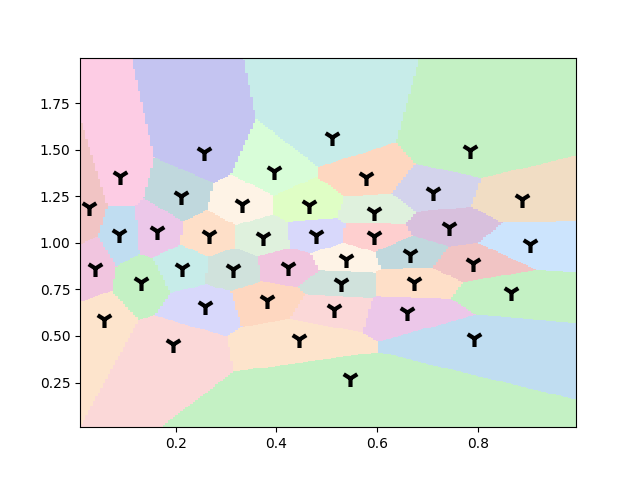}  
   \end{tabular}
    \caption[]{\em $P= {\cal N}\left(\begin{pmatrix}\tfrac 12 \\ 1 \end{pmatrix};\tfrac{1}{4}I_2\right)$. Bregman functions (from left to right and top to bottom):  $F(x)=|x|^2$, $F(x)=\sum_{i=1,2}\log(1+e^{a x_i})$ (SoftPlus, $a=1$), $F(x) =\sum_{i=1,2}e^{a x_i}$ (exponential, $a= 1$ and $a=-1$). With truncated $P$: 
$F(x) = -\sum_{i=1,2}\log x_i$ (Itakura--Sa\"ito), $F(x)=\sum_{i=1,2}x_i\log x_i$ (Kullback--Leibler),
\label{fig:2DN0I2}}
\end{figure}
\bibliographystyle{plain} 
\bibliography{these_artic.bib}

\end{document}